\documentclass[11pt]{article}
\usepackage{amssymb}
\usepackage{amsthm}
\usepackage{amsmath,amsfonts,amssymb}
\usepackage{diagbox}
\usepackage[mathscr]{eucal}
\usepackage{exscale} 
\usepackage{natbib} 
\usepackage{bm}
\usepackage{eqlist} 
\usepackage[dvipsnames]{color}
\usepackage[Lenny]{fncychap}
\usepackage{multirow}
\usepackage{graphicx}
\usepackage{algpseudocode}
\usepackage{algorithm}
\usepackage{caption}
\usepackage{hyperref}
\usepackage{framed}
\usepackage{color} 
\usepackage{booktabs}
\usepackage{longtable}

\usepackage{framed}
\usepackage{algpseudocode}

\usepackage{tabu} 
\usepackage{tikz}

\usepackage[utf8]{inputenc} 
\usepackage[english]{babel} 
\usepackage[T1]{fontenc}    
\usepackage{amssymb,amsmath,amsfonts} 
\usepackage{braket} 
\usepackage{caption}
\usepackage{algorithm}

\hypersetup{
	colorlinks=true,
	linkcolor=blue,
	citecolor=blue,
	citebordercolor={1 1 0},
}

\newtheorem{theorem}{Theorem}[section]

\newtheorem{definition}[theorem]{Definition}
\newtheorem{remark}[theorem]{Remark}

\font\bigbf=cmbx10 scaled \magstep3

\begin{document}

	\title{\bigbf Robust identification of drive-by sensing ride-hailing market with Points of Interest monitoring under uncertain ride demand} 
	\date{} 
	\author{Binzhou Yang \textsuperscript{a}$\thanks{Corresponding author, e-mail: ybz@my.swjtu.edu.cn ;}$  
	\quad Bin Shuai \textsuperscript{a} 
	\quad Rui Wang \textsuperscript{b}  
	\\\\
	\textit{\small a  School of Transportation and Logistics, Southwest Jiaotong University, China}\\
	\textit{\small b  College of Air Traffic Management, Civil Aviation Flight University of China, China}\\ 
	}
     
	\maketitle 
	
	\begin{abstract} 
		Taxi-based mobile sensing has emerged as a cost-efficient paradigm for large-scale urban environmental monitoring. In practice, both passive and active sensing strategies are adopted by ride-hailing platforms. Passive sensing, conducted during passenger-serving trips, is constrained by stochastic and spatially imbalanced ride demand, leading to limited and uneven coverage. Active sensing, executed by vacant taxis, provides greater control over sensing operations but incurs additional operational costs, and is thus typically treated as a supplementary strategy. To address the inefficiencies of the conventional ``passive-first, active-second'' paradigm, which may delay the monitoring of critical Points of Interest (POIs) and increase system-wide costs, we propose a Distributionally Robust Optimization (DRO)-based framework for active mobile sensing. First, we develop an enhanced A*-based drive-by sensing routing policy that integrates vehicle-task matching while capturing both global routing efficiency and local sensing opportunities. Second, we formulate the active sensing problem as a DRO model that explicitly accounts for uncertainty in ride demand through an ambiguity set of probability distributions, enabling robust decision-making for vacant taxi routing and vehicle-task matching. The proposed framework is evaluated on both static and dynamic mobile sensing settings using real-world data. Computational results demonstrate that our approach achieves better performance in terms of sensing coverage, operational cost, and robustness compared to benchmark strategies, highlighting the value of integrating distributional robustness into taxi-based sensing operations. 		
	\end{abstract}
	
	
	\section{Introduction}\label{secIntro}  
	By integrating Drive-by Sensing (DS) into urban transport vehicles, various sensing scenarios can be realized by leveraging the hosts' high mobility and low maintenance costs. These scenarios include air quality monitoring, traffic state estimation, noise monitoring, urban heat island effect analysis, parking availability detection, and built environment assessment \citep{MJKCXGT2010, PDB2012,  Gao2016, BAD2017, AD2017, Cruz2020a, Cruz2020b, TRMSO2020, GQDRJ2022}. Meanwhile, following a decade of rapid growth, the ride-hailing industry has seen a marked slowdown in the U.S. and China, its two largest markets, with fierce competition driving platforms to explore new value-added services beyond passenger transportation, such as providing food or freight delivery service, engaging with the third-party integrators \citep{WWWSGY2020,LGHFD2025}. In this context of DS, ride-hailing vehicles offer distinct advantages over fixed or dedicated sensing fleets, including dense urban coverage, long operating hours, and the potential for coordinated sensing task allocation.  
	 
	Due to platform dispatch protocols, monitoring is primarily limited to analyzing the spatio-temporal distribution of completed trips and the cruising of vacant taxis. While existing literature has explored sensor deployment, operational interventions, and incentive mechanisms for vehicle‑based sensing \citep{ZMLL2015, BBLARB2016, YHLL2025}, the majority of these works assume static or dedicated fleets and neglect the dynamic, demand‑driven trajectories of ride‑hailing vehicles. The resulting blind spots undermine the completeness of urban coverage data, introducing uncertainty and risk into subsequent urban management and governance \citep{BAD2017,OASSR2019,ADFS2021,LGHFD2025}. In contrast to prior work on sensing coverage for road area sources (1 km × 1 km grids in road networks), this study focuses on sensing coverage for point sources at monitoring blind spots furthermore, that is, Points of Interest (POIs) for road network monitoring. Thus, a fundamental challenge in existing research is the tension between leveraging opportunistic vehicle mobility for road segment monitoring and meeting the targeted spatio-temporal monitoring tasks, see Section~\ref{2.1} for further details. Indeed, on a ride-hailing platform with DS capabilities, aggregate POIs monitoring demand is often evaluated at a strategic level. Such an assessment, however, inherently overlooks the dynamics, uncertain conditions of real-time ride demand, which influences the POI monitoring.  
	
	To accomplish the sensing tasks throughout the monitoring time domain, the prevailing design of the DS ride-hailing market, which focuses on maximizing sensing gain via continuous order dispatch \citep{YHLL2025}, must be supplanted by a framework that incorporates uncertainties inherent to ride-hailing operational characteristics, see Section~\ref{2.2} for details. The realization of such a framework faces critical issues: (1) accurate matching and routing planning for occupied taxis, as generated by the platform's operational objectives for a mixed DS fleets; and (2) task assignment and monitoring routes for vacant instrumented taxis. To address the above two key issues, the vehicle‑passenger and vehicle-task matching schemes developed in this paper establish that its monitoring metric adopts priority‑based matching \citep{I2026,WJDL2026,KS2026}. The purpose of the former is to use passive sensing to maximize sensing utility while maintaining service levels (e.g., customer wait time and matching rate) for regular riders by prioritizing trip requests. On this basis, the latter aims to complete as many POI monitoring tasks as possible through active sensing, while achieving high sensing efficiency in road monitoring coverage over a short period. Under active sensing, the vacant instrumented taxis are controllable, the platform performs task assignment and route planning, and incurs the corresponding operation costs to incentivize POI monitoring. This paper subsequently proposes a robust framework account for vacant instrumented taxi routing and vehicle-task matching, which explicitly characterizes the associated uncertainties and POI monitoring cost. The technical approach of the proposed framework is summarized as follows. 
	
	Guided by the “passive-first, active-second” paradigm, in each decision-making epoch, trip requests are distributed to instrumented taxis first to maximize sensing utility, and then to non-instrumented taxis to meet ride demand. To supplement the above passive sensing, an A*-based drive-by sensing routing policy is implemented through active POI-monitoring. The matching and routing pairs of vacant instrumented taxis and tasks are constructed to characterize the interplay between local sensing opportunities and global routing efficiency. In detail, integrating the dual spatial scales of grid-level coverage and network-level routing, as many monitoring tasks as possible are assigned to vacant instrumented taxis in a time-efficient manner that maximizes local sensing utility. The final generated vehicle–task matching and routing pairs guarantee global efficiency. Subsequently, to address the uncertainties inherent in the DS ride-hailing platform's monitoring operations (e.g., daily travel demand fluctuations), a Distributionally Robust Optimization (DRO) model is employed and adapted to manage the monitoring costs associated with the designed routes for vacant instrumented taxis. Finally, a data-driven approach is employed to calibrate the uncertainty set within the DRO framework, enabling it to applicable for different operational control resolutions.
	
	This study offers the following major contributions:
	
	\textbf{Concept and modeling:} This work proposes a feasible enhanced A*-based drive-by sensing route service policy to determine the matches between vacant instrumented taxis and POI monitoring tasks, thereby ensuring local opportunistic sensing and global routing efficiency. Simulations confirm that the generated matching and routing pairs can complement the data collection blind spots while achieving high timeliness in road sensing coverage.
	
	\textbf{Theoretical properties:} This paper presents the first application of the DRO technique within the DS ride-hailing market. The DRO technique provides an approach for optimizing the operation cost incurred by POI monitoring with robust outcomes. It captures the interactive effects of the platform's uncertain ride demand on POI monitoring, without requiring assumptions about the individual or joint distributions of these uncertain variables.
	
	\textbf{Practical significance:} The proposed data-driven uncertainty set calibration is highly practical for the mixed-fleet DS ride-hailing market, enabling the DRO model to operate at different control resolutions (e.g., predetermined and adaptive). Furthermore, the calibration method can be readily generalized to account for varying data availability under day-to-day, within-day, and seasonal conditions.
	
	The remainder of this paper is structured as follows. Section \ref{secLR} reviews the relevant literature. Section \ref{secMe} details the proposed methodology, a simulation-based case study is presented in Section \ref{secNum}. Finally, Section \ref{secDis} discusses the framework's potential for enhancing the operational performance of a DS ride-hailing market with a mixed fleets, and Section \ref{secCon} provides concluding remarks.  
	
	\section{Related work}\label{secLR} 
	
	In this section, we provide a research survey of sensing operations for taxi fleets (Section \ref{2.1}), POI monitoring with uncertain ride demand (Section \ref{2.2}) and robust optimization approaches in transportation service (Section \ref{2.3}). Notably, under the uncertainty of POI sensing, our study addresses a research gap in the DS ride-hailing market with the aim of optimizing monitoring cost expenditures. 
	
	\subsection{Sensing operations for taxi fleets}\label{2.1} 
	
	Data collected from instrumented vehicles, with minimal interference to fleet operations, exemplifies the opportunistic sensing capabilities of taxi fleets. Utilizing empirical taxi data from Rome, \cite{BBLARB2016} showed that 80\% of the downtown area could be covered within 24 hours by a fleet of 120 vehicles. Highlighting the suitability of taxi fleets for scouting on-street parking availability, \cite{BAD2017} reported that about 500 taxis were sufficient to cover one-ninth of road segments in San Francisco, despite the spatio-temporal heterogeneous distribution of movement patterns. \cite{OASSR2019} examined how well taxi fleets can sense in several major cities by means of a ball-in-bin model, verifying considerably limited sensing capabilities due to the heterogeneous distribution of taxi trajectories. Regarding the optimal subset selection of drive-by sensing taxi fleets.\cite{ZYL2018} jointly optimized coverage across vehicles and per-task sensing reliability, framing vehicle selection as a bi-objective problem. \cite{YFLRL2021}, accounting for host vehicle speeds and their effect on sensing quality, designed a spatio-temporal, road-network-based vehicle selection model for traffic condition monitoring.
	  
	Operational interventions and incentives have been proposed to enhance the sensing coverage performance of taxi fleets. A primary focus is on routing interventions for given origin-destination pairs. For instance, \cite{Masutani2015} sought to maximize sensing quality by providing routing recommendations to vehicles. An $\varepsilon$-perturbed route set was proposed in \cite{ADFS2021}, generated through deviations from the A-star-computed shortest path, targeting optimized spatial coverage for a taxi fleet. Similarly, \cite{GQ2024} provided routing guidance to vacant taxis during their cruising phase. \cite{LGHFD2025} classified taxi types and utilized instrumented taxis for monitoring targeted objectives via the shortest path; however, their methodology fails to capture the utility of road segment sensing coverage. For incentive schemes, a scheduling-and-pricing model was proposed in \cite{FJLQG2021} to reward drivers whose routes stay within  a tolerable detour bound, thereby promoting coverage over broader areas. \cite{XCPJZN2020} formulated both routing advice and monetary incentives to maximize the sensing gain of vacant taxis cruising for their next customer. These incentive-based strategies were optimized using an efficient algorithm named iLOCuS. Furthermore, \cite{C2020} addressed the problem by matching riders with vacant taxis while providing routing suggestions, aiming to maximize sensing coverage under a limited incentive budget. This framework was integrated with a prediction-based actuation system called PAS. 
	
	Previous research has primarily focused on enhancing the overall sensing capabilities of taxi fleets, often overlooking specific target sensing objects. Consequently, the collected spatio-temporal data may lack statistical significance for targeted monitoring tasks. To address the demand for POI monitoring, this study leverages vacant instrumented taxis, which exhibit high controllability and flexibility. These vehicles are tasked to either complement the sensing capabilities of other fleets or to perform highly-targeted sensing missions. This objective is achieved by continuously coordinating sensing demand with travel demand through the implementation of an enhanced A* drive-by sensing route policy. 
	
	\subsection{POI monitoring with uncertain ride demand } \label{2.2}
	Instrumented taxis collect data on the urban physical environment by scanning road segments, whether engaged in  customer transport or in vacant cruising. The collected data serves as an input for governmental urban management. 
	
	Existing researches have focused on two main approaches to accomplish monitoring tasks: utilizing the inherent mobility of vehicles during transportation services, or executing predefined monitoring assignments. On leveraging the inherent mobility of vehicles, Despite demonstrating the suitability of taxi fleets for large-scale urban sensing (e.g., on-street parking availability as in \cite{BAD2017}), research also underscores a critical flaw: reliance on opportunistic, mobility-driven data collection leads to incomplete road coverage and spatio-temporal biases, which in turn create substantial uncertainty in the monitoring data. \cite{BAD2017} reported that approximately 500 taxis could cover one-ninth of road segments in San Francisco, yet this coverage remains inherently uneven due to heterogeneous movement patterns. Complementing this, \cite{OASSR2019} verified that the sensing capability is considerably limited precisely because of the uneven distribution of taxi trajectories. Together, these studies indicate that the very mobility which enables sensing also guarantees gaps in coverage, generating uncertainty that must be explicitly addressed. Regarding the predefined monitoring tasks, \cite{ADFS2021} compare the spatio-temporal coverage of probe vehicles routed by the $A^{*}$ algorithm against that achievable under the $RA^{*}_{\varepsilon}$ algorithm, across varying $\varepsilon$ values. This comparison is conducted on datasets of real-world taxi trajectories, considering specified source ($s$) and target ($t$) points. \cite{LGHFD2025} present the distribution pattern of monitoring POIs within the Manhattan network and employs a transportation-service-compatible shortest route planning strategy for vacant instrumented taxis.
	
	From the strategic level of the DS ride-hailing platform. A key limitation of prior work is its failure to account for the POI monitoring with uncertain ride demand when utilizing instrumented taxis for road segment sensing. This problem results in inherent gaps in coverage and persistent data blind spots. Consequently, the effectiveness of subsequent monitoring-data-based urban governance is compromised. To address this issue, this paper formulates and solves a priority-based matching and routing problem that jointly considers road segment sensing, ride-hailing services and POI monitoring demand with uncertainty.
	
	\subsection{Robust optimization approaches in transportation} \label{2.3} 
	
	Robust optimization is a methodology that addresses uncertainty by seeking solutions that remain feasible under the worst-case realization within a predefined uncertainty set. Since the publication of several seminal works in the late 1990s, robust optimization has seen substantial developments in both its theoretical foundations and practical applications. Compared to stochastic programming, robust optimization mitigates the computational intractability that can arise from the high dimensionality of uncertain parameters. A key advantage of robust optimization is that it relaxes the stochastic programming assumption that precise probability distributions of uncertainties are known, an assumption that often does not hold in practice. Furthermore, as decision-makers are often risk-averse rather than risk-neutral, the conservative nature of robust optimization makes it particularly applicable.
	
	DRO, first proposed by \cite{SH1957} in the context of a newsvendor problem, integrates concepts from both robust and stochastic optimization. The DRO framework optimizes the worst-case expected value of the objective function over an ambiguity set of probability distributions, which is characterized by certain known parameters of the uncertainty. A significant body of research is devoted to the construction of different ambiguity sets. A common approach involves defining the ambiguity set using fixed-moment information, such as the mean and variance. A comprehensive survey of this area is provided by \cite{GS2010}. Other methodologies include the work of \cite{CSS2007}, who incorporated directional deviations to refine the description of the distribution set, and \cite{DY2010}, who investigated a general DRO framework where the mean and covariance of the underlying uncertainties are subject to uncertainty. Both robust and distributionally robust optimization are widely applied in transportation, including route choice \citep{QSSY2016}, user equilibrium \citep{AMMN2015}, network design \citep{NL2016}, air traffic flow pattern identification \citep{SHMO2017}, routing optimization \citep{JQS2016,ZZLS2021}, facility location \citep{A2017}, and network sensor location \citep{ZFM2018}. 
	
	Existing research indicates that DRO techniques are capable of evaluating additional expenditures, such as monitoring costs, within a robust framework. However, a research gap remains in evaluating  the extra operational cost associated with monitoring POIs using vacant instrumented taxis. To bridge this gap, this paper constructs an uncertainty set for the extra monitoring cost under a designated monitoring pattern. By leveraging DRO techniques, we aim to minimize the expected cost expenditure under the worst-case scenario, thereby effectively capturing the potential uncertainties associated with the mixed DS ride-hailing platform operation under uncertain ride demand. Furthermore, this work provides an extensible framework for the ride-hailing market, with potential contributions to the strategic management of monitoring costs.
	
	\section{Method}\label{secMe}   
	Presented in this section are several elements fundamental to the later discussion of the robust identification method under uncertain ride demand, namely the background materials of DS ride-hailing market operation with POIs monitoring in Section \ref{3.1}, which includes the quantification of sensing externality in Section \ref{subsecDSquant}, a priority-based matching process in the context of DS in Section \ref{Mp} and the drive-by sensing coverage problem with POI monitoring in Section \ref{NG}, the enhanced A* drive-by sensing route service policy for POIs monitoring in Section \ref{3.2}, DRO extension of the deterministic formulation in Section \ref{DRO3.3}.

	Table \ref{tabmath} lists the notations to be used in the remainder of the paper.

	\setlength\LTleft{0pt}
	\setlength\LTright{0pt}
	\begin{longtable}{@{\extracolsep{\fill}} r p{12.0cm}}
	\caption{Mathematical notation}\label{tabmath} \\
	\hline
	\multicolumn{2}{l}{Sets and indices} \\
	\hline
	\endfirsthead
	
	\caption*{Table~\thetable{} (continued): Mathematical symbols} \\
	\hline
	\multicolumn{2}{l}{Sets and indices} \\
	\hline
	\endhead
	
	\hline
	\multicolumn{2}{r}{Continued on next page} \\
	\endfoot
	
	\hline
	\endlastfoot
	
	$\mathscr{A}$    & Driver set, whose elements are indexed by \(a\in \mathscr{A}\) \\
	$\mathscr{B}$    & Rider set, whose elements are indexed by \(b\in \mathscr{B}\) \\
	$\mathscr{A}^*$  & Matched-driver set\\
	$\mathscr{B}^*$  & Matched-rider set \\
	$T$              & Sensing-interval set, whose elements are indexed by \(t\in T\)\\
	$\mathcal{N}$    & Set of nodes \\
	$\mathcal{G}$    & Set of grids \\
	$\mathfrak{R}_G$ & Set of $G\text{-}\mathrm{route}$ \\
	$\mathfrak{R}_N$ & Set of $N\text{-}\mathrm{route}$ \\
	\hline
	\multicolumn{2}{l}{Parameter and constant definitions} \\ 
	\hline
	$l=(a,b)$        & A candidate pairing of driver \(a\) with rider \(b\) \\
	$\Gamma$          & Total sensing utility  \\
	$\Theta$     & Sensing utility  \\ 
	$\vartheta_{ab}$ & Sensing externality of a matched trip $l=(a, b)$ \\
	$\tau _{ab}$     & Pick‑up distance from the current location of driver $a$ to rider $b$ \\
	$w_{ab}$        & The customer wait time of a matched trip $l=(a, b)$ \\
	$N_{1}$          & Instrumented taxi fleets size \\
	$N_{0}$          & Non‑instrumented taxi fleets size \\
	\hline
	\multicolumn{2}{l}{Decision variables} \\
	\hline
	$x_{ab}$         & A binary indicator that is \(1\) when driver \(a\) and rider \(b\) are matched, and \(0\) otherwise \\
	
	\end{longtable}

	\subsection{Background of DS ride-hailing market operation with POIs monitoring}\label{3.1}
	\subsubsection{Quantifying drive-by sensing (DS) externality}\label{subsecDSquant}
	For a discretized spatio-temporal sensing framework, let $G$ denote the set of unit grid cells (each $1\,\text{km} \times 1\,\text{km}$) partitioning the target region, and let $T$ index discrete sensing intervals (e.g., hourly). The overall sensing utility is: 
	\begin{equation}
		\Gamma = \sum_{t \in T} \mu_t \sum_{g \in G} s_g \, \Theta_{g,t}(\Lambda_{g,t}),
		\label{eqnPhi}
	\end{equation}
	with normalized temporal weights $\mu_t$ ($\sum_{t \in T} \mu_t = 1$) and spatial weights $s_g$ ($\sum_{g \in G} s_g = 1$) reflecting sensing priorities.  For grid $g \in G$ during interval $t \in T$, the definition of the sensing utility is: 
	\begin{equation}
		\Theta_{g,t}(\Lambda_{g,t}) = (\Lambda_{g,t})^{\rho}, \quad \rho \in (0,1),
	\end{equation}
	where $\Lambda_{g,t}$ gives the count of distinct vehicles visiting grid $g$ during interval $t$. The exponent $\rho$ captures diminishing marginal returns, discouraging overconcentration of visits. We set $\rho = 0.2$ for our air-quality monitoring case study \citep{JHL2023, HJNLL}. 
	
	In online driver-rider matching, each sensing interval $t$ is subdivided into decision epochs $s_i$ ($i = 1, \dots, n$). For a candidate match $l = (a, b)$ considered at epoch $s_i$, its marginal utility is:
	\begin{equation}
		\vartheta_{ab} = \sum_{t \in T} \mu_t \sum_{g \in G_{ab}} s_g \left( (\Lambda_{g,t,s_i} + 1)^{\rho} - (\Lambda_{g,t,s_i})^{\rho} \right),
		\label{zetadr}
	\end{equation}
	where $G_{ab} \subset G$ denotes the set of grids passed through by trip $l$, and $\Lambda_{g,t,s_i}$ denotes the cumulative vehicle count for grid $g$ in interval $t$ as of epoch $s_i$.  
	
	\subsubsection{A priority-based matching process in the context of DS}\label{Mp}  
	
\begin{figure}[h]
		\centering
		\includegraphics[scale=0.50]{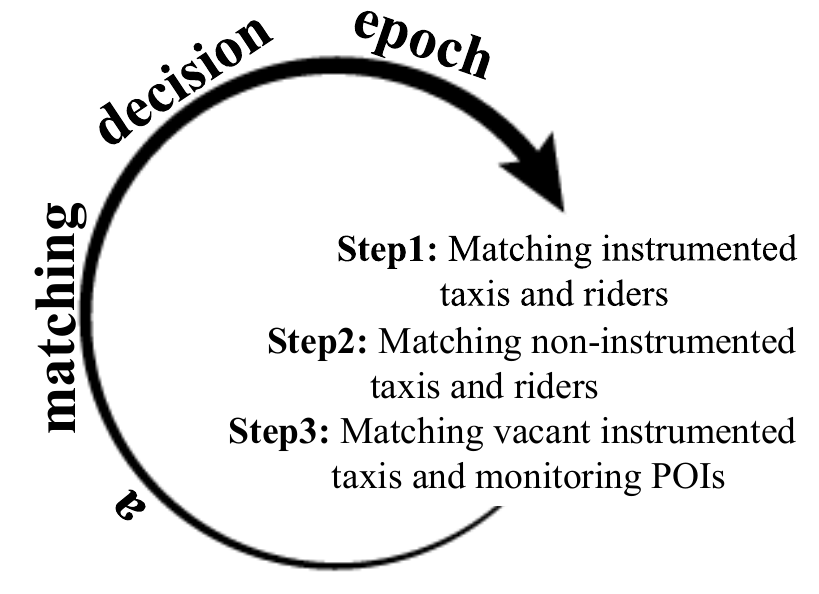}\\
		\caption{Matching process in a decision epoch}\label{figmodelsplit} 
 		\end{figure}

	Taxis are divided into instrumented (fixed size \(N_1\)) and non-instrumented (fixed size \(N_0\)) fleets. In the context of DS, the matching process is priority-based to serve the platform’s  operations \citep{I2026,WJDL2026,KS2026}. As illustrated in Figure \ref{figmodelsplit}, orders are first assigned to instrumented taxis to maximize sensing gain, then to non‑instrumented taxis to optimize operational metrics with the objective of minimizing customer wait time. Any remaining vacant instrumented taxis are dispatched to monitoring POIs to fulfill sensing tasks. 
	
	In Step 1, the objective focuses on the sensing gain $\vartheta_{ab}$ for a matched trip, which is calculated using Eq. \eqref{zetadr}. This calculation assumes that the route is known in advance and follows the shortest path. The matching model targets the maximization of the total sensing externality, formulated as:

	\begin{equation}\label{mopt1}
		\hbox{[Sensing utility]}~~   \max_{x=(x_{ab})_{a\in\mathscr{A},b\in\mathscr{B}}} S=\sum_{b\in\mathscr{B}} \sum_{a\in\mathscr{A}}  \vartheta_{ab} x_{ab}
	\end{equation}
	\begin{equation}\label{mopt2}
		\hbox{s.t.}~\sum_{b\in \mathscr{B}}x_{ab}\leq 1,\quad \forall a\in\mathscr{A} 
	\end{equation} 
	\begin{equation}\label{mopt3}
		\sum_{a\in \mathscr{A}}x_{ab}\leq 1,\quad \forall b\in\mathscr{B}
	\end{equation}
	\begin{equation}\label{mopt4}
		x_{ab} \in \{0, 1 \} , \quad \forall a\in\mathscr{A},\,b\in\mathscr{B}
	\end{equation} 
	Here, $x_{ab}$ denotes a binary decision variable, equal to 1 if driver $a$ and rider $b$ are matched,
	
	In Step 2, the objective shifts to minimizing the customer wait time $w_{ab}$ for a matched trip $l=(a,b)$, assuming travel along the shortest route. Consequently, the matching model for this step is formulated as follows:
	
	\begin{equation}\label{mopt5}
		\hbox{[Customer wait time]}~~   \min_{x=(x_{ab})_{a\in\mathscr{A},b\in\mathcal{B}}} W=\sum_{b\in\mathscr{B}} \sum_{a\in\mathscr{A}}  w_{ab} x_{ab}
	\end{equation}
	\begin{equation}\label{mopt6}
		\hbox{s.t.}~\sum_{b\in \mathscr{B}}x_{ab}\leq 1,\quad \forall a\in\mathscr{A} 
	\end{equation} 
	\begin{equation}\label{mopt7}
		\sum_{a\in \mathscr{A}}x_{ab}\leq 1,\quad \forall b\in\mathscr{B}
	\end{equation}
	\begin{equation}\label{mopt8}
		x_{ab} \in \{0, 1 \} , \quad \forall a\in\mathscr{A},\,b\in\mathscr{B}
	\end{equation} 
	The binary variable \(x_{ab}\) indicates whether driver \(a\) and rider \(b\) are matched, taking the value 1 if so. The continuous variable $w_{ab}$ denotes the resulting customer wait time for this match.	
	
	In Step 3, the matching process involves vacant instrumented taxis and monitoring POIs, with the specific objectives and routing strategies detailed in Section \ref{NG}.
	
	\subsubsection{The drive-by sensing coverage problem with POI monitoring}\label{NG} 
	The  Drive-by Sensing Coverage (DSC) problem is proposed by  \cite{HJNLL} firstly. Unlike classical facility location with fixed sensors, DSC employs instrumented taxis as mobile sensors, extending the traditional framework \citep{RHS2021,WWWZQ2021,VRCT2024}. The problem of routing instrumented taxis for data collection has two key aspects: (i) sensing utility is accrued via grid intersections rather than node/arc visits, and (ii) utility follows a diminishing‑returns function $\Theta_{g,t}(\cdot)$. The primary challenge in coordinating grid-level coverage with network-level routing lies in their coupling across disparate spatial resolution scales.
	
	\begin{figure}[h]
		\centering
		\includegraphics[scale=0.26]{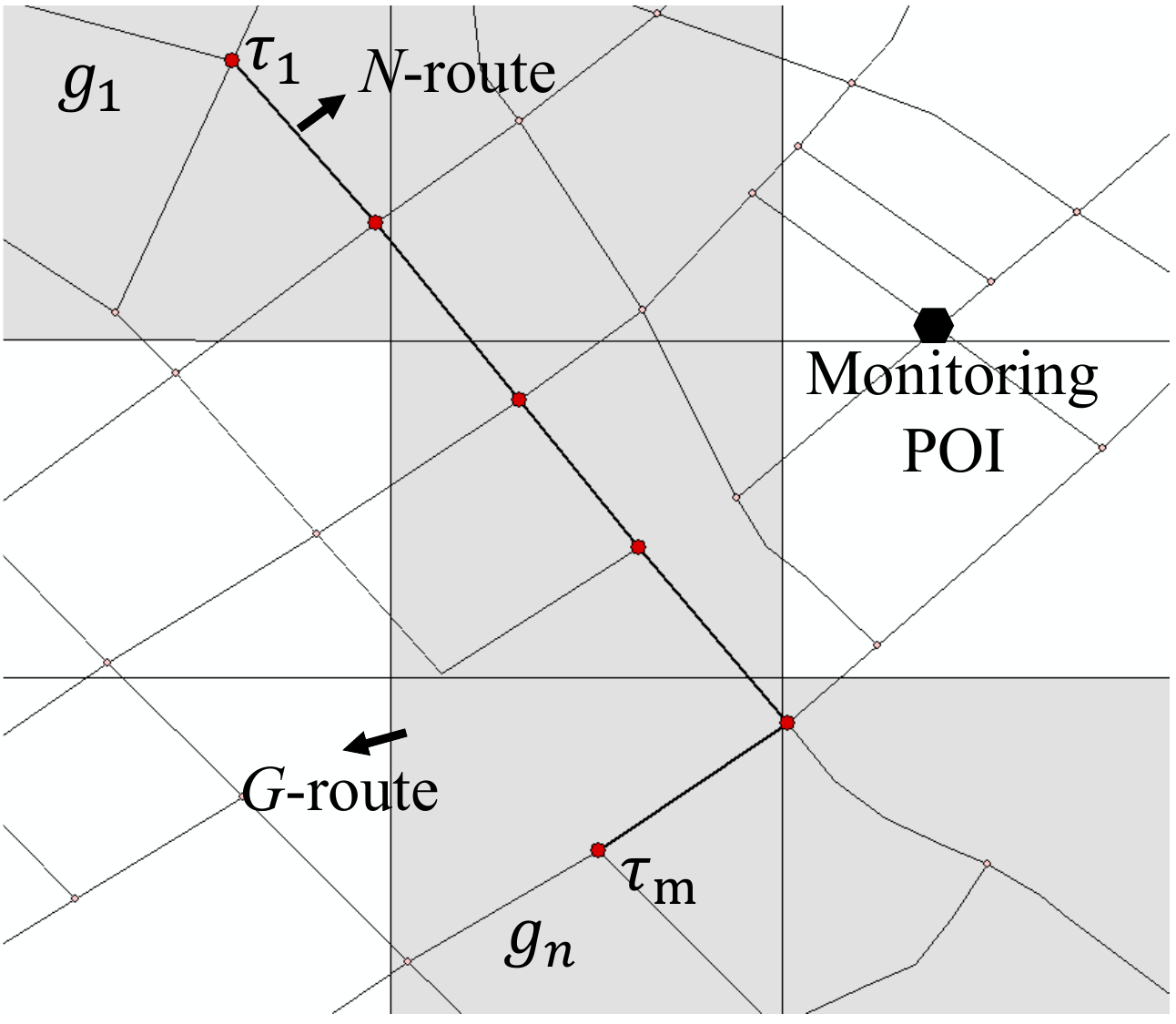}
		\caption{Drive-by sensing coverage with POI monitoring. Adapted from ‘‘Exploring the sensing power of mixed vehicle fleets'' by \cite{HJNLL}, Transportation Research Part B: Methodological. 190, 103066 (https://doi.org/10.1016/j.trb.2024.103066). Copyright 2024 by Elsevier.}
		\label{figNGroute}
	\end{figure}
	
	In this paper, considering monitoring blind spot coverage to supplement monitoring data, we extend the DSC monitoring problem to the DSC problem with POI monitoring. As illustrated in Figure \ref{figNGroute}, a $N\text{-}\mathrm{route}$ (where $N\text{}\mathrm{}$ stands for network) refers to a conventional route represented as a sequence of nodes $\kappa_N = (\tau_1, \dots, \tau_m) \subset \mathcal N$. A $G\text{-}\mathrm{route}$ (where $G\text{}\mathrm{}$  stands for grid) refers to a conceptual route expressed as the sequence of grids traversed, $\kappa_G = (g_1, \dots, g_n) \subset \mathcal G$  \citep{HJNLL}. A node $\tau \in \mathcal N$ is said to belong to a grid $g \in \mathcal G$ (denoted $\tau \in g$) if it lies within that grid, and Figure \ref{figNGroute} shows a POI to be monitored through a marked node. Therefore, accurately characterizing the spatio-temporal patterns of POI monitoring demand is essential for effectively utilizing the available fleet of vacant instrumented taxis. The design of efficient matching and routing algorithms for POI monitoring at dual spatial scales is detailed in Section \ref{3.2}.
	 
	\begin{remark}
		Building on the deterministic formulation presented 
		in Section \ref{3.2}, which assumes prior and accurate knowledge of POI monitoring demand, we provide the monitoring path policy for POIs. From an operational perspective, the generated monitoring path of vacant instrumented taxis for POIs directly determines the POI monitoring cost (i.e. operating cost of vacant instrumented taxis for POI monitoring). Nevertheless, the resulting solutions of monitoring cost may lack robustness against potential variations in demand distribution caused by multiple sources of uncertainty. 
		
		To overcome this limitation, we adopt a data‑driven distributionally robust optimization (DRO) approach. The detailed formulation is provided in Section~\ref{DRO3.3}. DRO addresses uncertainty in both parameters and their underlying distributions, which are difficult to estimate accurately from limited or interdependent data. Unlike conventional robust optimization, it avoids excessive conservatism by optimizing against the worst-case distribution within a plausibly specified family. We propose a data-driven \(\phi\)-divergence-based calibration method to construct the uncertainty set. Leveraging empirical statistics and risk preferences with sufficient data, the proposed test defines distributions that are close to the empirical estimate.
	\end{remark}

	\subsection{Enhanced A* drive-by sensing route service policy for POI monitoring} \label{3.2}  
	
	\subsubsection{POI monitoring routing for a vacant instrumented taxi}

	For a single instrumented taxi, we consider a matching and routing problem defined from an origin node \( \tau \) in the road network. The overall objective is to determine a destination node \( \tau_m \) (which represents a monitoring POI) along with the corresponding $N\text{-}\mathrm{route}$  and  $G\text{-}\mathrm{route}$ to the target grid, so as to gather sensing utilities in a time‑efficient manner. Given the origin node \( \tau_1 \in g_1 \in \mathcal G \) in the road network and taking into account Eq. (\ref{zetadr}), the problem amounts to identifying vehicle-task matching pairs with the utility time‑efficiency function \( \varphi(\cdot) \).	 	 
	
	\begin{equation} \label{ueffi}    
		\varphi(a\tau)= \frac{ {\textstyle \vartheta_{\bar{\kappa}(\tau_{1};\kappa_{G})}} }{T_{\bar{\kappa}(\tau_{1};\kappa_{G})}} 
	\end{equation}

	\noindent where the $N\text{-}\mathrm{route}$ fully implements the $G\text{-}\mathrm{route}$  \( \kappa_G = (g_1, \dots, g_n) \in \mathfrak{R}_G \) with starting node \( \tau_1 \in g_1 \), which is denoted \( \bar{\kappa}(\tau_1;\kappa_G) \in \mathfrak{R}_N \). In each matching decision epoch, monitoring POIs are assigned to instrumented taxis through the following matching model.  
	
	\begin{equation}\label{ueffiopt1}
		\max_{x=(x_{a\tau})_{a\in\mathscr{A},\tau\in\mathcal{N}}}Q=\sum_{\tau\in\mathcal{N}} \sum_{a\in\mathscr{A}} \varphi(a\tau) x_{a\tau} 
	\end{equation}
	\begin{equation}\label{ueffiopt2}
		\hbox{s.t.}~\sum_{\tau\in\mathcal{N}}x_{a\tau}\leq 1,~\forall a\in \mathscr{A}
	\end{equation} 
	\begin{equation}\label{ueffiopt3}
		\sum_{a\in \mathscr{A}}x_{a\tau}  \leq 1,~\forall \tau\in\mathcal{N}
	\end{equation}
	\begin{equation}\label{ueffiopt4} 
		x_{a\tau} \in\left \{ 0 ,1\right \}, ~\forall a\in\mathscr{A},\,\tau\in\mathcal{N}  
	\end{equation}
	\noindent Here, the binary decision variable \(x_{a\tau}\) indicates whether driver \(a\) is matched with monitoring POI \(\tau\), with the corresponding $N\text{-}\mathrm{route}$    \(\bar{\kappa}(\tau_{1};\kappa_{G})\) generated by a modified A* algorithm. The utility time‑efficiency \(\varphi(a\tau)\) is defined as the ratio of the weighted marginal sensing utilities accumulated along the traversed $G\text{-}\mathrm{route}$  \(G_{\bar{\kappa}(\tau_{1};\kappa_{G})}\) to the total travel time \(T_{\bar{\kappa}(\tau_{1};\kappa_{G})}\). Maximizing \(\varphi(a\tau)\) incentivizes taxis to prioritize grids with high marginal sensing utility, thereby accumulating greater utility within a shorter time window. In other words, the proposed algorithm is designed to incentivize vacant instrumented taxis to complete POI monitoring by preferentially traversing grids with high marginal sensing utility.	This procedure is formally summarized in Table \ref{tabmath1}.

	\subsubsection{POI monitoring routing for multiple vacant instrumented taxis}
	To account for overlapping search areas, we adopt a sequential optimization procedure that iteratively updates the marginal utilities of traversed grids until convergence or a maximum iteration limit is reached. Routes are planned within each discrete time window using the window-specific total sensing utility \(\Gamma\), thereby avoiding redundant cross-window coverage and reducing the computational burden of long‑horizon planning. 
	
	Algorithm 2’s routing plan in Table \ref{tabmath2} is extended over the full time horizon \(T\) by converting the routes \(R_{1}^{A},\ldots,R_{n^{A}}^{A}\) into repeatable round trips, with each one-way leg executed within a discrete time window \(t\). Unlike conventional vehicle routing problems that assume fixed depots, instrumented taxi routes may originate from any road network node, provided that monitoring POI demand is satisfied (see Section \ref{4.5} for an example).
	
	\subsection{DRO extension of the deterministic formulation}\label{DRO3.3}
	
	\subsubsection{Introduction to the general data-driven DRO framework}\label{framework}   
	Let $\textbf{Q}= \left({\textbf{Q}_{t}}: t\in T \right)$ 
	denote a random vector that affects the POI monitoring during each discrete sensing interval $t\in T$. In a more general context, $\textbf{Q}$ can represent any uncertain variables or vectors related to ride demand uncertainty, either as exogenous inputs (e.g. city size, weather and population density) or endogenous variables (e.g. trip fare, fleet availability). 
	
	For each vacant instrumented taxi, the generated POI monitoring route (detailed in Section \ref{3.2}) is used to compute the POI monitoring cost, based on a unit-distance compensation of 1.5 CNY/km provided to the driver. The total monitoring cost borne by the mixed DS ride-hailing platform is the sum of the monitoring costs associated with all vacant instrumented taxis.  Given a fixed variable $x$, which denotes the fleet size of instrumented taxis ($N_{1}$) in this study, the monitoring cost is considered as the objective function or Key Performance Indicator (KPI), denoted $KPI(x,\textbf{Q})$, which is to be minimized and treated as a random variable parameterized by $x$. 
	  
	\begin{equation}
		KPI(x, \mathbf{Q}) = 1.5 \times \sum_{i=1}^{N_1} L_i(\mathbf{Q})
	\end{equation}
	\noindent where $x = N_{1}$  and $L_i(\mathbf{Q})$ denotes the per-unit-time POI monitoring route length (km) for instrumented taxi $i$. 
	
	This randomness arises from the stochasticity of $\textbf{Q}$, and is assumed to follow an unknown probability distribution $D^{*}(x)$. In the absence of exact knowledge of $D^{*}(x)$, we instead consider an uncertainty set $\mho{(x)}$ comprising candidate probability distributions that are consistent with the available data via the   $\phi$-divergence test. Hence, given that the KPI (i.e., the objective function) is to be minimized, the DRO formulation is presented as follows:
	
	\begin{equation}\label{objective}
		\underset{x}{\text{min}}\underset{D\in \mho{(x)} }{\text{max}} E_{D}[KPI(x,\textbf{Q})]
	\end{equation}	
	where $E_{D}$ denotes the expectation KPI under the candidate probability distribution $D(x)$. The DRO model is formulated as a min-max optimization problem, which is non-convex and infinite-dimensional. 
	
	To ensure the practical solvability of the problem, we adopt a metaheuristic solution approach supported by a finite approximation of the PDF. Specifically, for the outer minimization problem in Eq. (\ref{objective}), metaheuristic search algorithms, such as the genetic algorithm (GA) and particle swarm optimization (PSO), can be employed to search for the optimal solution  $x^{*}$, relying solely on the objective function values (i.e., the results of the inner maximization problem). Thus, it suffices to determine the solution to the inner maximization problem for any given $x$. The detailed implementation procedure is presented as follows.	
	
	We consider a set of $M$ historical samples of the random vector \textbf{Q}, denoted  $({\textbf{Q}^{1}},{\textbf{Q}^{2}},...,{\textbf{Q}^{M}})$. Given the instrumented taxi fleet size $x$, the monitoring cost (i.e., the KPI) is computed based on the unit-distance compensation rule for vacant instrumented drivers (1.5 CNY/km). This computation is performed via simulation or a deterministic procedure, specifically the enhanced $A^{*}$ drive-by sensing route service policy proposed in this study.  Consequently, a sequence of KPI samples is obtained, denoted as { $KPI(x,\textbf{Q}^{1}), KPI(x,\textbf{Q}^{2}), ..., KPI(x,\textbf{Q}^{M})$ }, where each sample represents the monitoring cost derived from the $k$-th historical data sample $\textbf{Q}^{k}$ (with $k = 1, 2, ..., M$). This sequence is defined as:
	\begin{equation}
		f^{k} = KPI(x, \textbf{Q}^{k}), \quad 1 \le k \le M
	\end{equation}
	where each $f^k$ corresponds to a realization of the random variable $KPI(x, \textbf{Q}^{k})$. Without loss of generality, we assume these KPI samples are sorted in non-decreasing order, i.e., 
	$f^{1}\le f^{2}\le...\le f^{M}$.
	\begin{definition} 
		$\phi$-divergence is a distance measure between two vectors $\textbf{p}=(p_{1},p_{2},...,p_{m})\ge 0$ and $\textbf{q}=(q_{1},q_{2},...,q_{m})\ge0$ in $\textbf{R}^{m}$ and can be defined as:
		\begin{equation}\label{xiequation}
			I_{\phi}(p,q)=\sum_{j=1}^{m}q_{j}\phi(\frac{p_{j}}{q_{j}})
		\end{equation}
		
		If $p$ and $q$ are probability vectors, there exists constraints $\sum_{j=1}^{m}p_{j}=1$ and $\sum_{j=1}^{m}q_{j}=1$. $\phi(u)$ is a convex function for $u>0$.
	\end{definition} 
	
\begin{figure}[h]
		\centering
		\includegraphics[scale=0.5]{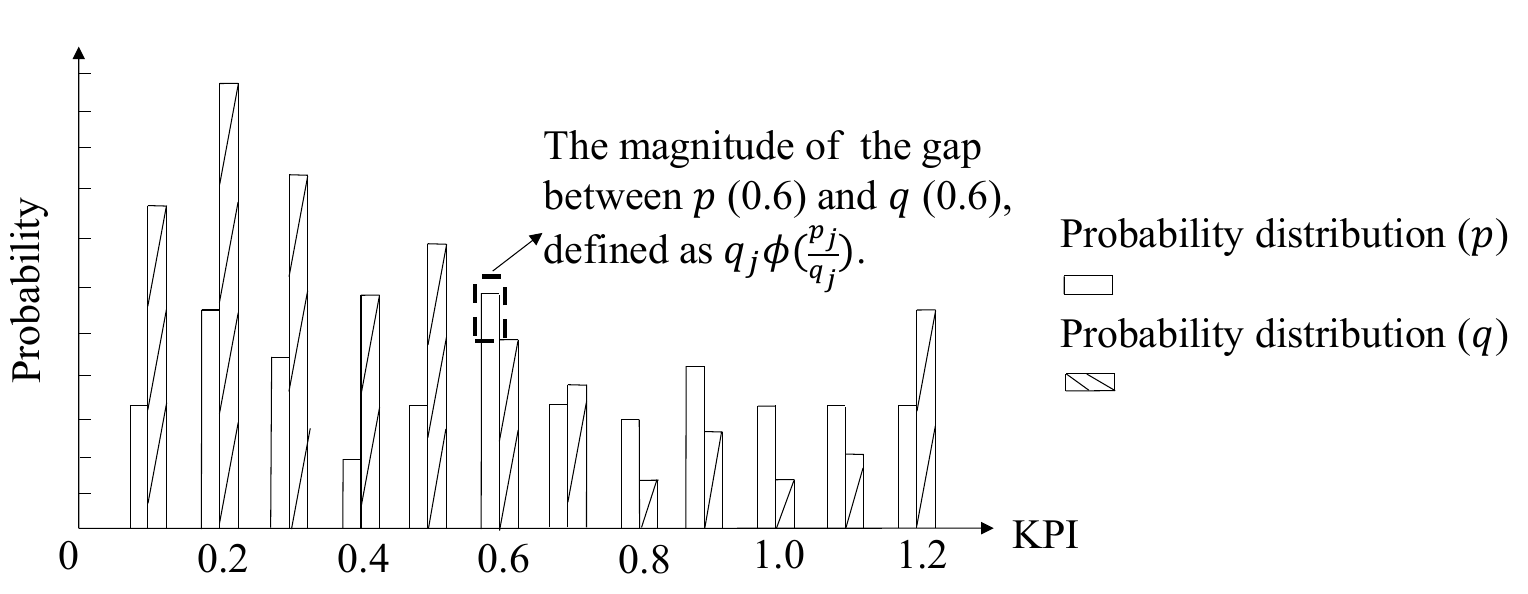}\\
		\caption{An illustration of phi-divergence.}\label{figofphidivergence} 
	\end{figure}

	Intuitively, $\phi$-divergence is a mathematical measure used to quantify the discrepancy between two PDFs. To facilitate understanding, two discrete probability distributions are illustrated in Figure \ref{figofphidivergence}. For a given scenario $j$ in Eq. (\ref{xiequation}), there exists a discrepancy between the probabilities associated with the two distributions, denoted $p_{j}$ and $q_{j}$, respectively. The magnitude of this discrepancy is defined as $q_{j}\phi(\frac{p_{j}}{q_{j}})$, and the $\phi$-divergence between the probability vectors $\textbf{p}$ and $\textbf{q}$ is defined as the sum of these terms over all scenarios $j=1,2,3,...,m$, i.e. $\sum_{j=1}^{m}q_{j}\phi(\frac{p_{j}}{q_{j}})$.
	
	Adopting a discrete probability framework, in contrast to a continuous one, reduces model complexity while retaining essential probabilistic characteristics. The corresponding uncertainty set is constructed as follows:
	
	\begin{equation}\label{Judgment condition} 
		\mho{(x)} =\left \{p_{j}=(p_{1},p_{2},...,p_{m})|I_{\phi}(p,q)\le \rho, \sum_{j=1}^{m}p_{j}=1,p_{j} \ge 0  \right \} 
	\end{equation}
	where $\rho$ is a predefined parameter representing the uncertainty level. The formulation of $\mho{(x)}$ can accommodate various $\phi$-divergence functions, allowing flexibility to align with decision-maker preferences. A comprehensive summary of these functions is provided in Table \ref{divergencecategory}.	
	
	The $\phi$-divergence is grounded in the following hypothesis testing framework:
	
	\begin{equation}
		H_{0}:D(x)=D^{*}(x), \quad H_{A}: D(x) \ne D^{*}(x) 
	\end{equation}
	where $D^{*}(x)$ is the unknown true distribution, and $D(x)$ is any nominal distribution.
	The null hypothesis $H_{0}$ is rejected at significance level $\alpha$ if the following holds \citep{BDDMR2013}:
	\begin{equation}
		T_{\phi}^{M}(p,q)=\frac{2M}{\phi^{''}(1)}I_{\phi}(p,q)
	\end{equation} 
	where $T_{\phi}^{M}(p,q) $ $\sim$ $ \chi_{m-1}^{2}$.  
	The rejection of the null hypothesis $H_{0}$ is jointly determined by the sample data and a prescribed quantile of a specific distribution. For a comprehensive derivation, see \cite{BDDMR2013}. The rejection rule is
	\begin{equation}
		T_\phi^M(p,q) > \chi_{m-1,\,1-\alpha}^2 
	\end{equation}
	with \(\chi_{m-1,\,1-\alpha}^2\) denoting the \((1-\alpha)\)-quantile of the chi-square distribution with \(m-1\) degrees of freedom.
	Consequently, the parameter \(\rho\) is defined as
	\begin{equation}
		\rho = \frac{\phi''(1)}{2M} \chi_{m-1,\,1-\alpha}^2,
	\end{equation}
	which depends on the sample size \(M\), the number of grid points \(m\), and the significance level \(\alpha\).
	
	Building upon this characterization of candidate distributions, we now focus on approximating the worst-case expectation, which corresponds to the inner problem of Eq.~(\ref{objective}). To this end, we discretize the probability density function (PDF) of the KPI. Let $L, U \in \mathbb{R}$ with $L < U$ be the lower and upper bounds of the KPI over the feasible set $\mathcal{Q}$ of the uncertain parameter $\mathbf{Q}$:
	\begin{equation}
		L = \inf_{\mathbf{Q} \in \mathcal{Q}} KPI(x, \mathbf{Q}), \qquad 
		U = \sup_{\mathbf{Q} \in \mathcal{Q}} KPI(x, \mathbf{Q}).
	\end{equation}
	
\begin{figure}[h]
		\centering
		\includegraphics[scale=0.49]{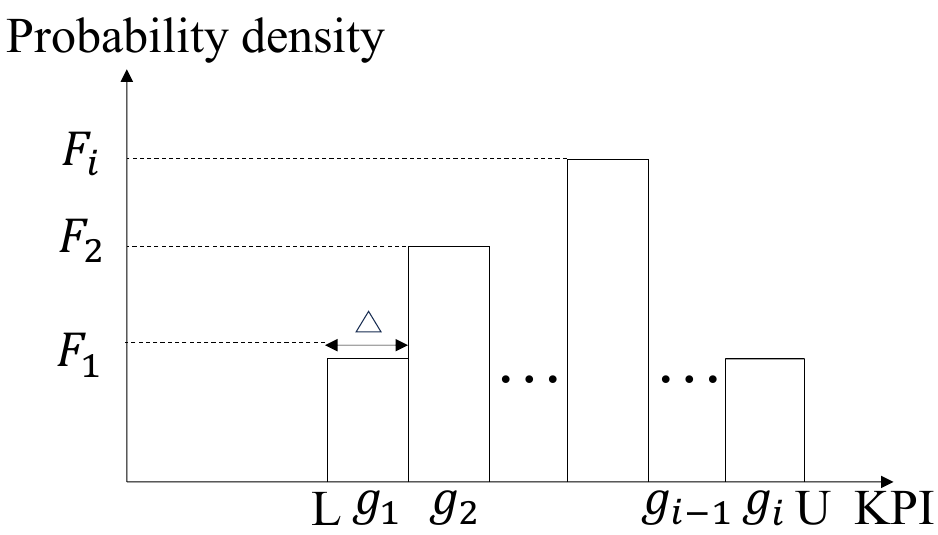}\\
		\caption{Discrete approximation of the PDF of KPI.}\label{figdiscretekpi} 
	\end{figure}

	The interval $[L, U]$ is uniformly partitioned into $w$ subintervals of length $\Delta = (U-L)/w$. The grid points $(g_i)_{1 \le i \le w}$ are placed at the centers of these subintervals, i.e., $g_i = L + (i-\frac12)\Delta$, as illustrated in Figure~\ref{figdiscretekpi}. The continuous PDF $F$ is then approximated by the discrete vector $\hat{\mathbf{F}} = (\hat{F}_i)_{1 \le i \le w}$, where $\hat{F}_i \approx F(g_i)$.
	
	The probability mass at grid point $g_i$ is given by $p_i = \hat{F}_i \Delta$, 
	which approximates the integral of the PDF over the $i$-th subinterval 
	using the rectangle rule (midpoint rule). These probability masses satisfy $\sum_{i=1}^w p_i = 1$ and $p_i \ge 0$. The uncertainty set $\mho(x)$ imposes additional constraints on the shape of the PDF, which translate into box constraints on the probability masses:
	\begin{equation}
		\underline{p}_i \le p_i \le \overline{p}_i, \quad \forall i,
	\end{equation}
	where $\underline{p}_i$ and $\overline{p}_i$ are the minimum and maximum feasible values of $p_i$ derived from $\mho(x)$.
	
	Consequently, the inner problem in Eq.~(\ref{objective}) can be approximated by the following linear program (LP):
	\begin{equation}
		\max_{D \in \mho(x)} \mathbb{E}_D[KPI(x, \mathbf{Q})] \approx
		\max_{\mathbf{p}} \left\{ \sum_{i=1}^{w} g_i p_i \;:\;
		\sum_{i=1}^{w} p_i = 1,\;
		p_i \ge 0,\;
		\underline{p}_i \le p_i \le \overline{p}_i,\; \forall i \right\}
		\label{eq:LP}
	\end{equation}
	
	\subsubsection{Selection of $\phi$-divergence function}
	The selection of an appropriate $\phi$-divergence function is a critical step, which should be guided by the characteristics of the available data \citep{BL2015, ZFM2018}. These adjustments are characterized by two classification rules: suppression and popping.
	\begin{definition}[Suppression]
		The ability to assign zero or reduced probability to scenarios that are present in the empirical sample but are deemed unreliable.
	\end{definition}
	
	\begin{definition}[Popping]
		The ability to assign positive probability to scenarios that are \emph{not} observed in the empirical sample.
	\end{definition}
	
	Specific instances of these divergence functions are summarized in Table \ref{divergencecategory}, with explanatory notes provided subsequently. It is important to note that this classification is not unique; alternative frameworks could be applied depending on the modeling context.
	
	Whether suppression is permissible depends on the quality of the historical data. Suppression Allowed: When prior data are of poor quality, e.g., collected over short durations, from uncensored road segments, or under atypical conditions such as severe weather or accidents, suppression prevents the model from overfitting to sparse or unrepresentative samples. The optimization then concentrates on more trustworthy observations. Suppression Disallowed: If the historical data are abundant and of high quality, all observed scenarios should be retained to avoid discarding valuable empirical information. Disallowing suppression yields a conservative model that fully exploits the available data.
	
	Whether popping is beneficial depends on the availability of qualitative information beyond the quantitative sample. Popping Allowed: Popping is desirable when qualitative sources, e.g., expert opinions, surveys, or domain knowledge, suggest the existence of plausible scenarios not captured by quantitative data (such as GPS traces or sensor readings). Enlarging the support of the distribution expands the family of candidate probability models and mitigates the risk of underestimating uncertainty. Popping Disallowed: If no theoretical or qualitative basis supports the existence of unobserved scenarios, or if decision makers intend to adhere strictly to the empirical distribution, popping should be disallowed. The support is then confined to the observed sample.
	
	The Kullback-Leibler divergence in Table \ref{divergencecategory} is selected in this paper because it rigorously rejects distributions exhibiting ‘‘popping'' (assigning probability to unobserved values). This ensures conservative, empirically supported models that control tail risks and prevent over‑optimistic extrapolation.
	
	\subsubsection{DRO under predetermined controls}\label{DROpc}
	At the pre-tactical level, the DRO framework informs fleet‑size decisions prior to daily operations, using demand forecasts as input. Constructing this formulation first requires an initial uncertainty set, which is empirically built from \(M\) days of historical trip requests following the procedure in Section~\ref{framework}.
	
	The monitoring cost differential, the extra operational cost incurred when instrumented vehicles are rerouted from passenger-seeking to sensing tasks, is treated as a random variable influenced by latent, interrelated factors regarding demand fluctuations. Instead of assuming a parametric distribution, we adopt a nonparametric, data‑driven representation using the bivariate structure:
	\begin{equation}
		\mathbf{\Xi}^{k} = \big( \xi^{k}_{t}, \zeta^{k}_{t} \big)_{t \in T}, \quad k = 1, \dots, M,
	\end{equation}
	where \(\xi^{k}_{t}\) denotes the baseline monitoring cost during interval \(t\) on day \(k\), and \(\zeta^{k}_{t}\) is an auxiliary cost modifier that depends on the operational context  \citep{YZH2026} \footnote{In this paper, the auxiliary cost modifier aligns with the subsidy incentive cost item, which is designed to incentivize instrumented drivers to continuously accept orders.}. The collection 
	\[
	\mathfrak{D} = \{\mathbf{\Xi}^{1}, \mathbf{\Xi}^{2}, \dots, \mathbf{\Xi}^{M}\}
	\]
	consists of \(M\) i.i.d. realizations of an unspecified random vector \(\mathbf{Q}\), providing the empirical basis for the subsequent distributionally robust treatment.
	
	To evaluate the objective in Eq.~\eqref{objective} for a candidate instrumented fleet size \(N_1\), we first retrieve historical trip records covering the designated sensing horizon. The enhanced A* routing protocol (Section~\ref{3.2}) is then applied to simulate per‑interval monitoring costs, yielding \(m\) independent KPI realizations, each representing the total monitoring cost incurred by the mixed‑fleet ride‑hailing system. A linear programming routine subsequently computes the worst‑case expectation from this empirical KPI sample, which serves as the robustness‑adjusted performance metric. Since the set of feasible fleet sizes is finite, the optimal configuration is identified by direct enumeration, for higher‑dimensional decision spaces, metaheuristic methods (e.g., genetic algorithms, simulated annealing) are employed to maintain computational tractability.
	
	\begin{remark}
		A key methodological feature of this pipeline is that it omits explicit behavioral or structural assumptions about DS ride-hailing market dynamics (e.g., trip generation). Instead, this information is incorporated via two complementary mechanisms: an empirically derived uncertainty set, and trajectory-level realism from the A* routing heuristic. This approach is analogous to Monte Carlo simulation. Rather than exhaustively modeling stochastic primitives, we simulate system realizations directly from historical data. The resulting empirical distribution of the objective is then embedded into the distributionally robust optimization (DRO) framework without parametric intermediaries.
	\end{remark}
	
	\subsubsection{DRO under adaptive controls based on rolling horizon}\label{DROad} 
	The DRO-based pre-tactical planning model in Section~\ref{DROpc} produces decisions from day-ahead demand forecasts that remain static throughout the operating day, offering no adaptability to real-time updates. In practice, however, demand predictions are continually refined at shorter look-ahead times (LAT), e.g., hourly, yielding rolling forecasts that are more accurate and inherently capture real-time factors, such as weather, that affect ride demand. To harness these progressively improved estimates, we adopt a rolling horizon framework that reapplies the DRO methodology from Section~\ref{DROpc} at each LAT-aligned interval with minimal modifications, allowing each decision to incorporate the latest forecast. The LAT itself is not fixed, it embodies a context-dependent trade-off between forecast precision and computational cost (e.g., DRO solution time) and can therefore be tailored to diverse operational environments.
	
	\section{Numerical test}\label{secNum}		
	This section presents a case study conducted in the Longquanyi District of Chengdu. It begins with a description of the trip demand configuration, simulation initialization, distribution of monitoring POIs, taxi operation modes, and the optimization method, followed by implementation details in Section \ref{4.1}. Subsequently, Section \ref{4.2} presents the results of the deterministic, DRO, and DRO with rolling horizon frameworks. Section \ref{4.3} shows the simulation results of key operational performance indicators. Finally, the performance of monitoring timeliness for POIs and the enhanced A* drive-by sensing routes are discussed in Sections \ref{4.4} and \ref{4.5}, respectively.
	
\begin{figure}[h]
		\centering
		\includegraphics[scale=0.450]{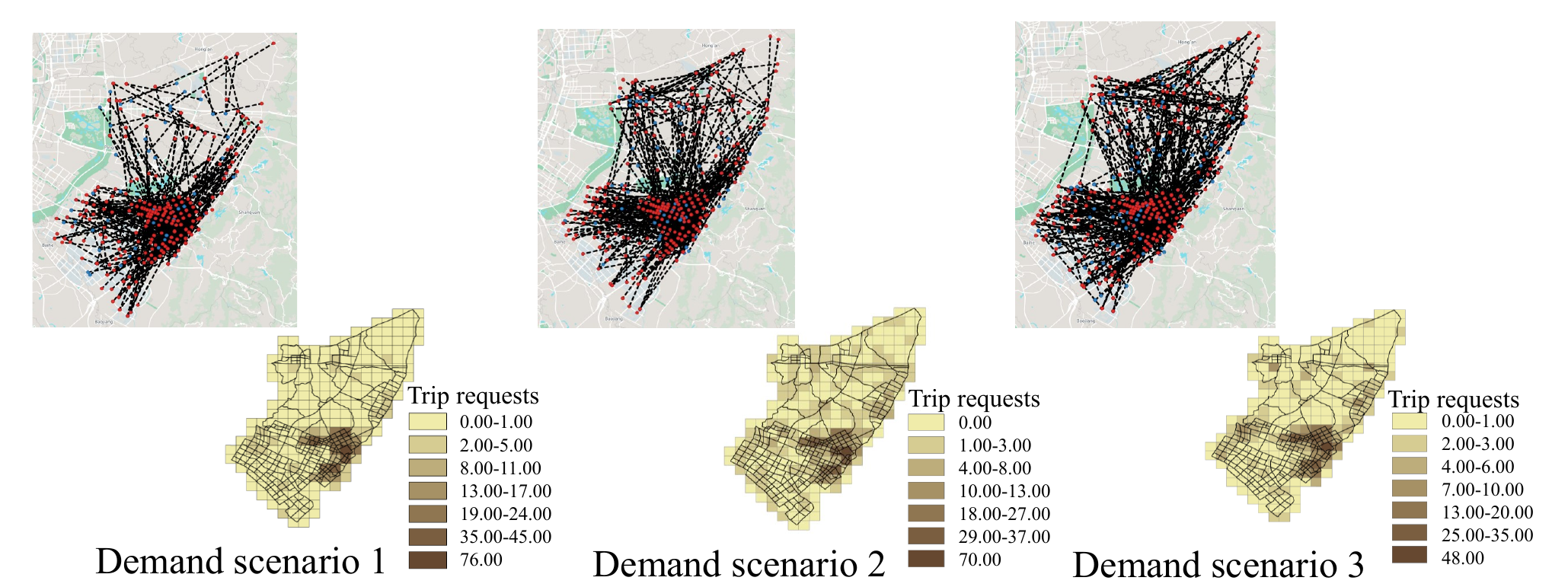}\\
		\caption{These three demand scenarios differ in the growing proportion of remote orders. The origin is marked in red, and the destination is marked in blue.}\label{figlmhdemand} 
	\end{figure}

	Trip Demand Configuration.  
	Trip demand is characterized using historical completed trip data collected between 8:00 and 12:00 from August 1 to December 31, 2021.  To capture the influence of demand patterns on system performance, multiple demand scenarios are constructed. Demand Scenario 1 serves as the baseline and is derived directly from actual completed trips. Two additional scenarios (Demand Scenarios 2 and 3) are synthetically generated by progressively increasing the proportion of remote trips, defined as trips ending in areas with low local demand, while maintaining the overall trip request volume at roughly 144 per hour. Thus, the three scenarios differ exclusively in the spatial distribution of trip demand, as illustrated in Figure \ref{figlmhdemand}. Induced demand and unmet demand are not considered in this study.	
	
	Simulation Initialization. 
	The simulation lasts four hours and is divided into four sensing intervals, one hour apiece. Each interval is evenly divided into 18 decision epochs. At the start of each epoch, the ride-hailing platform collects real-time trip requests and the locations of vacant instrumented taxis. 
	Every simulation run begins with all taxis randomly placed in the road network. Vacant taxis initially move toward areas of greater trip request density before being assigned rides.
	
	\begin{figure}[h]
		\centering
		\includegraphics[scale=0.34]{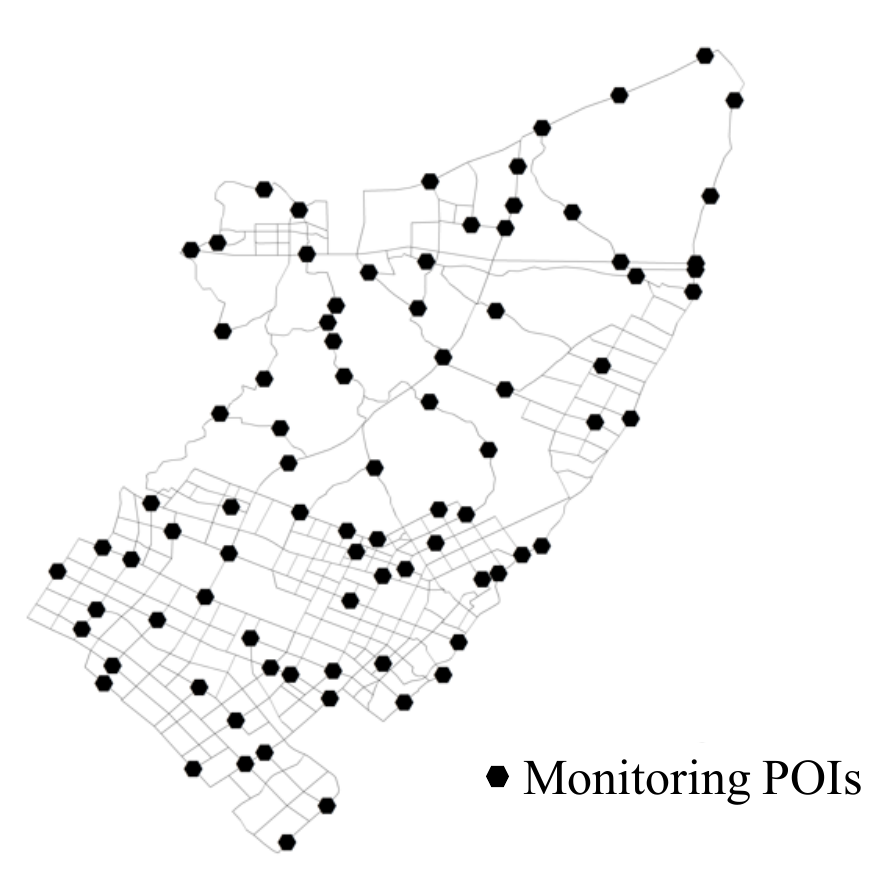}\\
		\caption{Monitoring POIs in Longquanyi District.}\label{monitoring POI} 
	\end{figure}
	
	Distribution of Monitoring Points. A total of 85 monitoring points of interest (POIs) are randomly selected as reference locations, as depicted in Figure \ref{monitoring POI}.
	
	Taxi Operation Mode.   
	The operating speed is set to a constant 35 km/h. At each decision epoch, the platform executes the matching mechanism described in Section \ref{Mp}, generating matched pairs and optimal routes. The system state is propagated to the subsequent epoch to ensure continuity in the dynamic simulation. Vacant taxis continue cruising when no ride assignments are received.
	
	Optimization Methods.  
	Depending on data availability and the scope of optimization methods, three variants of the proposed framework are implemented. The deterministic approach applies pre-determined control strategies without accounting for uncertainty. The DRO-based approach (Section \ref{DROpc}) also adopts pre-determined control but explicitly considers demand uncertainty. The DRO\(_{RH}\) approach (Section \ref{DROad}) incorporates adaptive control with a rolling horizon scheme, wherein a one-hour demand prediction is provided at the beginning of each time window. All methods are evaluated using the 85 randomly selected monitoring POIs.
	
	\subsection{Implementation details} \label{4.1}	
	The proposed matching procedure, executed within each decision epoch, is implemented using the following techniques and specifications (see Figure \ref{figmodelsplit} for details):
	\begin{enumerate} 
		\item Matching Radius: A search radius of 2 km is applied for driver-rider matching. Consequently, driver-rider pairs separated by more than 2 km are excluded from the matching process.  
		
		\item Step 1 (Sensing-Optimized Matching): In Step 1 of Figure \ref{figmodelsplit}, only the service leg of a matched trip (excluding the pick-up phase) contributes to the sensing externality objective in \eqref{mopt1}. This design is meant to prevent excessively long pick-up distances resulting from pursuing sensing utilities. Under the constraint of maximizing the sensing externality, each rider is matched to the nearest available driver to minimize customer wait time.
		
		\item Step 2 (Non-Instrumented Taxi Matching):  We maximize the number of matches between non-instrumented taxis and riders. This step thereby improves the overall matching rate after Step 1.
		
		\item Step 3 (POI Monitoring Task Matching): With the maximum number of matches between vacant instrumented taxis and monitoring POIs as the primary objective to fulfill more monitoring tasks, we then maximize utility time-efficiency. Note that if a POI was already included in the route of a pair matched in Step 1, it is considered a completed monitoring task and is excluded from Step 3.
	\end{enumerate} 
	
	\subsection{Results of Deterministic, DRO, and DRO with rolling horizon frameworks}\label{4.2} 		
	Table \ref{monitoringcostresults} summarizes the training and testing performance of the Deterministic, DRO, and DRO$_{RH}$ models in terms of monitoring cost. The Design Key Performance Indicator (Design KPI) presented in this table is defined as the average monitoring cost per instrumented taxi per unit time, which incorporates both the operational cost of the enhanced A* drive-by sensing route service policy and the subsidy incentives provided to the instrumented taxis. 
	
	The ``Optimal size'' column lists the optimal fleet size (in vehicles) identified by each method. The ``Design KPI'' column reports the actual performance metric when the optimal fleet size is applied to the realized demand data, as opposed to the predicted demand used for optimization.

	\begin{table}[!htp] 
		\caption{Training and testing results of the three optimization frameworks}
		\label{monitoringcostresults}  
		\resizebox{1.\textwidth}{.9in}{ 
			\begin{tabular}{ccccccccccc}
				\hline 
				\multirow{3}{*}{\begin{tabular}[c]{@{}c@{}}Testing\\ time\end{tabular}}                                           & \multirow{3}{*}{\begin{tabular}[c]{@{}c@{}} Demand\\ distribution\end{tabular}} & \multirow{3}{*}{\begin{tabular}[c]{@{}c@{}}Optimization\\ method\end{tabular}} & \multicolumn{6}{c}{Instrumented taxis fleet size} & \multirow{3}{*}{\begin{tabular}[c]{@{}c@{}}Optimal\\ size\\(vehicles)\end{tabular}} & \multicolumn{1}{c}{\multirow{2}{*}{\begin{tabular}[c]{@{}c@{}}Design KPI \\ (CNY/h/vehicle) for\\ the optimal size\end{tabular}}} \\ \cline{4-9}
				\\
				&       &   &10   & 20    & 30    & 40   & 50   & 60    & \multicolumn{1}{c}{}                                                                                            \\  
				\hline
				\multirow{9}{*}{\begin{tabular}[c]{@{}c@{}}8:00-12:00,\\ collected\\ from 1 \\ Aug to 31\\ Dec, 2021\end{tabular}} & \multirow{3}{*}{\begin{tabular}[c]{@{}c@{}} Demand \\ Scenario 1\\ \end{tabular}}   & \begin{tabular}[c]{@{}c@{}}Deterministic\end{tabular}                            
				& 50.935 &  23.999  & 15.341 & 11.633  & 8.374  & 6.963  & 60  & 6.963                                                                                                         \\ \cline{3-11} 
				&      & DRO  & 36.053 &18.932 &12.590 & 9.227 & 6.815 &4.759 & 60 &4.759                                                                                                           \\ \cline{3-11} 
				&      & \begin{tabular}[c]{@{}c@{}}DRO$_{RH}$  \end{tabular}                               
				& 35.717  & 18.077  & 11.585  & 8.321  & 6.281  & 5.137  & 60  &      5.137                                            \\ \cline{2-11}  
				& \multirow{3}{*}{\begin{tabular}[c]{@{}c@{}} Demand \\ Scenario 2  \end{tabular}} & \begin{tabular}[c]{@{}c@{}}Deterministic\end{tabular}  &  45.186 & 25.083 &14.388 & 10.452 &8.207 & 6.715& 60 & 6.715                                                    \\ \cline{3-11} 
				&       & DRO       &   33.629 & 18.914 &11.691& 8.146 & 6.623 & 5.271 & 60 & 5.271                           \\ \cline{3-11}   
				&         & \begin{tabular}[c]{@{}c@{}}DRO$_{RH}$\end{tabular}            &32.416 & 18.330 & 10.924 & 7.433 & 5.917 & 4.854 & 60 & 4.854                                                                                   \\ \cline{2-11} 
				& \multirow{3}{*}{\begin{tabular}[c]{@{}c@{}} Demand \\ Scenario 3  \end{tabular}}   & \begin{tabular}[c]{@{}c@{}}Deterministic\end{tabular}       &40.938 & 26.171 & 16.216 & 11.147&8.252 & 6.057 & 60 & 6.057                                                   \\ \cline{3-11}  
				&      & DRO        & 33.153 & 18.603 & 12.167 & 8.847 & 6.348 &4.995 & 60 & 4.995                     \\ \cline{3-11} 
				&      & \begin{tabular}[c]{@{}c@{}}DRO$_{RH}$\end{tabular}          &  31.967 & 17.343 &11.142 &7.644& 5.744 & 4.297 & 60 & 4.297                                  \\ \hline
		\end{tabular}}
	\end{table}
	
	Specifically, across three demand scenarios with fleet sizes ranging from 10 to 60 vehicles, DRO$_{RH}$ achieves superior performance in 17 out of 18 comparisons (94.4\%), with an average KPI reduction of 0.3–1.3 CNY/h/vehicle compared to the standard DRO approach. The advantage is most pronounced under high-demand conditions (Demand  Scenario 3), where DRO$_{RH}$ consistently outperforms DRO across all fleet sizes with KPI improvements of up to 1.26 CNY/h/vehicle at 20 vehicles (18.603 vs. 17.343 CNY/h/vehicle).  Notably, the only exception occurs in the lowest demand scenario (Demand Scenario 1) with the largest fleet size (60 vehicles), where DRO slightly outperforms DRO$_{RH}$ by a marginal 0.378 CNY/h/vehicle (4.759 vs. 5.137 CNY/h/vehicle)—a mere 7.4\% difference.  This single counterexample actually reinforces the method's validity, as it aligns with theoretical expectations: when the fleet size is sufficiently large relative to demand, temporal demand variations are smoothed out, making the pooled estimation approach (DRO) marginally more efficient. These results demonstrate that DRO$_{RH}$ effectively captures time-varying demand patterns under realistic resource constraints, validating its robustness for practical implementation.
	\begin{figure}[h]
		\centering
		\includegraphics[scale=0.35]{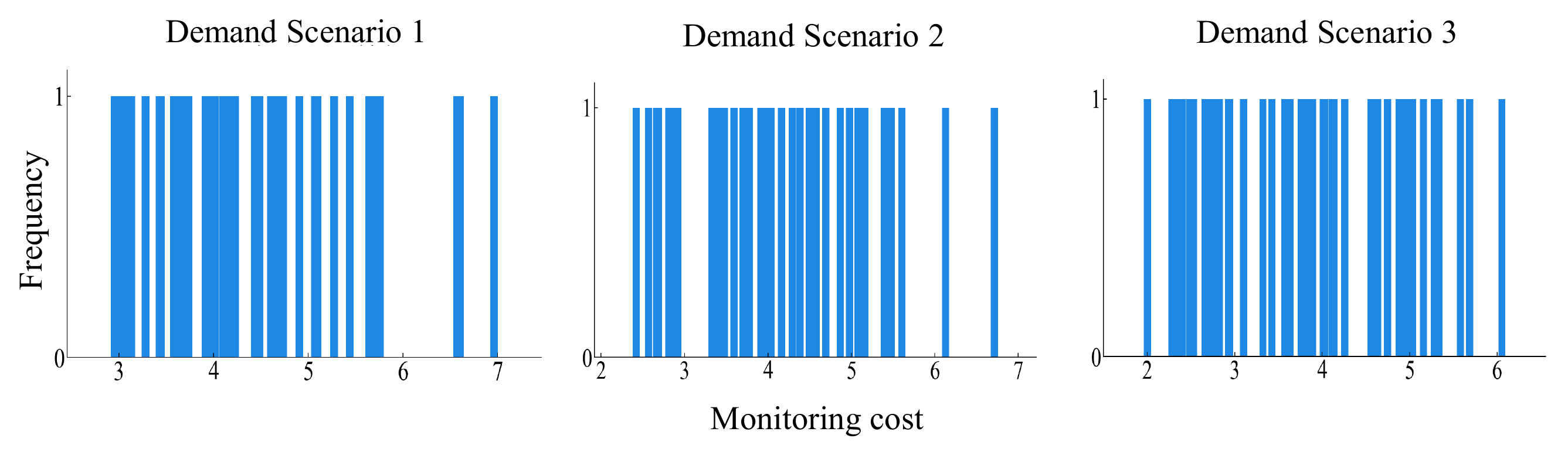}\\
		\caption{Histograms of the ``Design KPI'' under the optimization method of DRO$_{RH}$ approach, based on results from 10 simulation runs across three demand scenarios. The X-axis represents the per-unit-time monitoring cost per instrumented taxi (CNY/h/vehicle), and the Y-axis shows the corresponding frequency (count) in the simulation tests.}\label{figsd} 
	\end{figure}
	
	The results indicate that the DRO$_{RH}$ model outperforms the other two approaches. It achieves this by dynamically adjusting the fleet size (between 10 and 60 vehicles) in response to updated travel demand over time. This evaluation demonstrates the effectiveness of the proposed controls under real-world uncertainty. The DRO approach mitigates the impact of outliers in a tunable manner through the selection of $\alpha$. This adjustability enables a less conservative interpretation of uncertainties, which in turn leads to superior system performance, as evidenced in Figure \ref{figsd}.

	\subsection{Simulation results on operational performance KPIs}\label{4.3}

	The evaluation of operational metrics is organized into three aspects. Level of service comprises two indicators: (1) \textit{Matching rate} (\%) and (2) \textit{Average wait time} (min). The former refers to the proportion of riders who are assigned to drivers, while the latter denotes the average pick-up time for matched pairs. Sensing externality includes two indicators: (3) \textit{Sensing utility} \(\Gamma\), as defined in Eq. \eqref{eqnPhi}, and (4) \textit{Grid coverage rate} (\%), defined as the mean percentage of grids covered at least once per hour. Monitoring capability is measured by (5) \textit{monitoring completion rate} (\%), defined as the overall completion ratio of monitoring POIs.

	Figure \ref{figsim3casesm} through Figure \ref{figsim3casess} present the operational performance across three Demand Scenarios for fleet sizes ranging from 10 to 60 vehicles.  Based on 10 independent simulation runs, the following observations are drawn from the figures:
	
	\begin{enumerate} 
		\item  With matching rate and average wait time as its measures, the level of service responds to changes in the instrumented fleet size. Demand Scenario 3 shows the lowest matching rate, as drivers tend to cruise in areas with concentrated order density. As shown in Figure \ref{figmodelsplit}, the multi-objective optimization of fixed prioritization strategies leads to rigid resource allocation. In a dynamic environment with arriving trip requests and monitoring POIs, the increase in instrumented taxis leads to periodic conflicts between the sensing objective and the customer-matching objective, resulting in non-monotonic oscillatory variations in the matching rate.

		\item  Demand Scenario 3, with the highest proportion of remote orders, also exhibits the highest average wait time. The average wait time shows a non-monotonic trend, initially increasing and then decreasing with fleet size, which stems from the optimization mechanism illustrated in Figure \ref{figmodelsplit}. With fewer instrumented taxis, the platform prioritizes sensing maximization, leading to longer customer wait times. As the fleet expands beyond a certain threshold, sufficient instrumented taxis are assigned under this objective to serve a significant share of travel demand. This increase in instrumented taxis alleviates the conflict between sensing and service objectives, allowing the subsequent customer wait time minimization module to assign nearer taxis more effectively, thereby reducing the overall average wait time.

		\item  The sensing externality (encompassing sensing utility and grid coverage rate) and monitoring capacity (indicated by monitoring completion rate) both improve as the size of instrumented taxi fleets increases. Demand Scenario 3 exhibits higher sensing utility and grid coverage rate, attributable to its higher proportion of remote orders.
	\end{enumerate} 	
	
	\subsection{Monitoring timeliness of POIs}\label{4.4}  
	
\begin{figure}[h]
		\centering
		\includegraphics[scale=0.65]{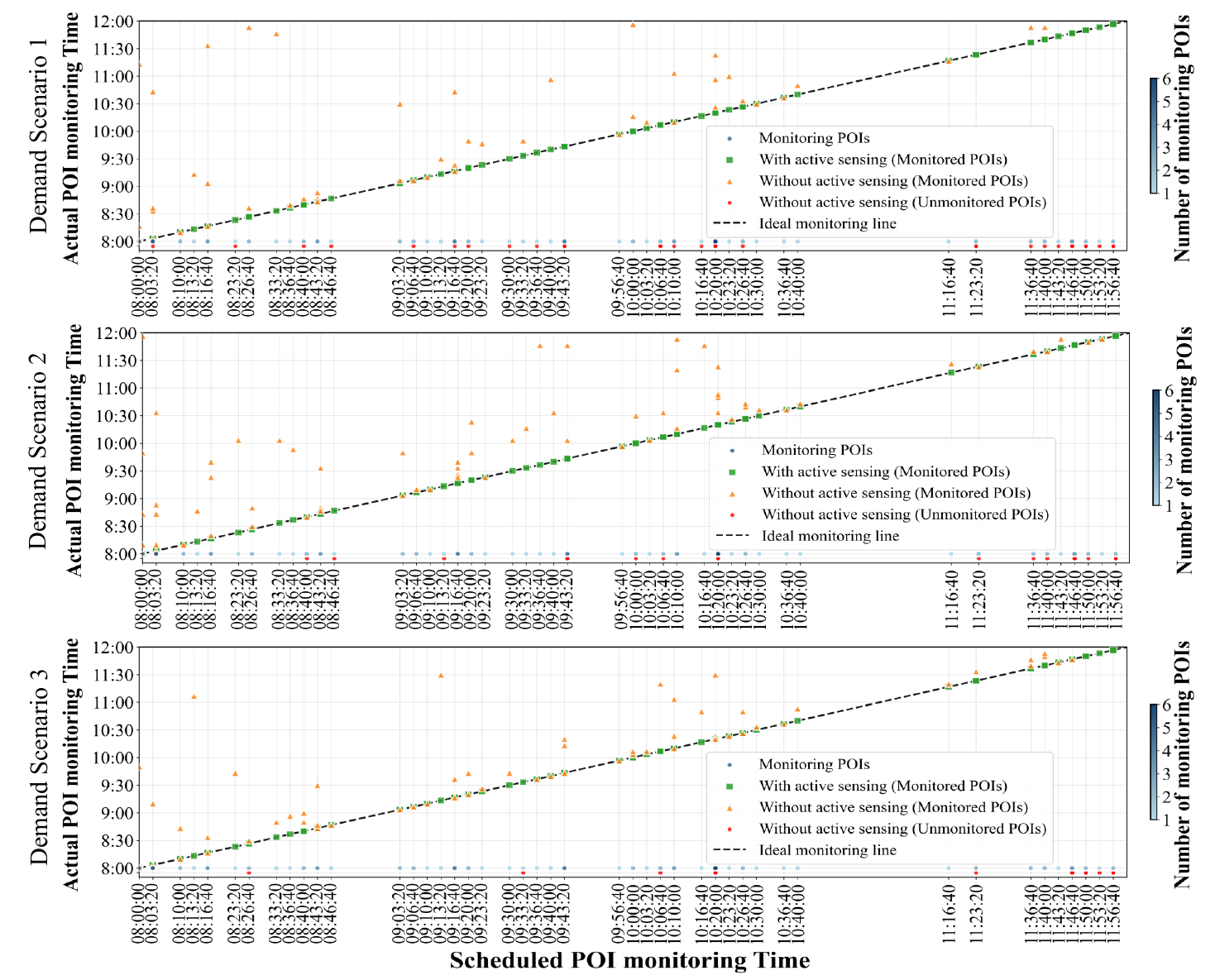}\\
		\caption{\textcolor{black}{Monitoring timeliness of POIs under Demand Scenario 1, 2 \& 3 (with regards to 60 instrumented taxis).}}\label{figsensingtime} 
	    \end{figure}

	Using the dispatch cycle from 8:00 to 12:00 as a reference, we statistically analyze the scheduled monitoring times of 85 POIs, as shown on the horizontal axis of Figure \ref{figsensingtime}. To verify that deploying vacant instrumented taxis for POI monitoring can achieve timely and effective monitoring coverage and data collection in blind spots, we examine the monitoring timeliness of these POIs under two modes: with active sensing and without active sensing. 
	
	The following conclusions are drawn:
	\begin{itemize}
		\item The blue dots represent the scheduled monitoring times of POIs. The darker the color, the greater the number of associated POIs. The red dots represent POIs that are not monitored within the calculated time domain. This situation occurs only in the ‘‘without active sensing'' mode.
		\item The diagonal line, where the x-coordinate equals the y-coordinate for every point on it, is defined as the ideal monitoring line. If a point defined by the actual monitoring time (vertical  axis) and the scheduled monitoring time (horizontal axis) falls exactly on this line, it indicates strong timeliness for that POI. As can be seen from Figure \ref{figsensingtime}, the actual POI monitoring times in the ‘‘without active sensing’’ mode are scattered more in the upper-left area above the ideal monitoring line rather than on the line itself. The timeliness of the ‘‘with active sensing’’ mode is better than that of the ‘‘without active sensing’’ mode.	
		\item When the order distribution is more balanced, there are fewer unmonitored POIs in the ‘‘without active sensing’’ mode. This is because POIs widely distributed across the measurement area (as shown in Figure \ref{monitoring POI}) are more likely to be covered and monitored by fleets of instrumented taxis that are occupied (i.e., serving passengers).
	\end{itemize} 

	\subsection{Enhanced A* drive-by sensing routes}\label{4.5}  	 
\begin{figure}[h]
		\centering
		\includegraphics[scale=0.37]{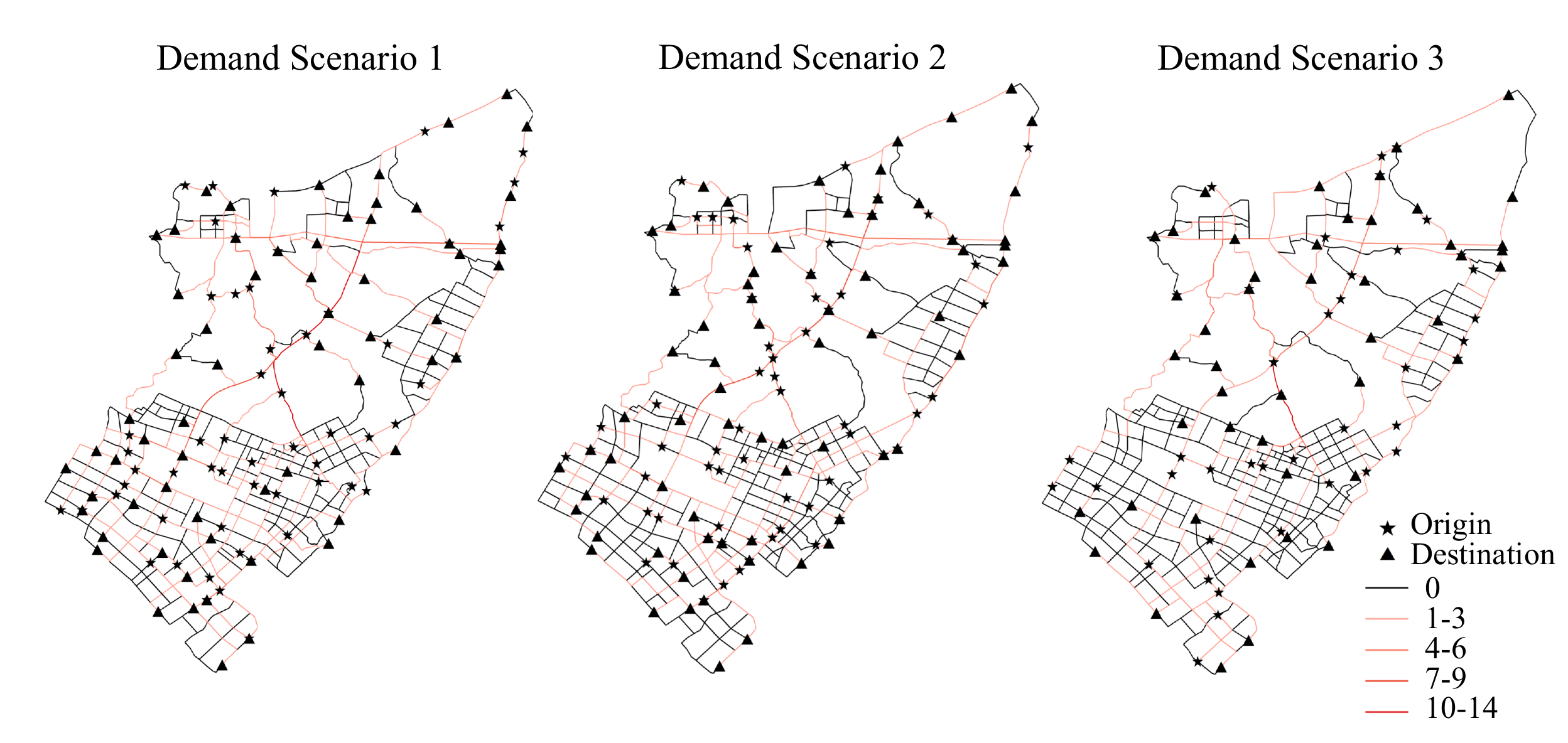}\\
		\caption{\textcolor{black}{Enhanced A* drive-by sensing routes under Demand Scenario 1, 2 \& 3 (with regards to 60 instrumented taxis).}}\label{figsensingroutes} 
	\end{figure}

	With reference to Table \ref{monitoringcostresults}, a fleet size of 60 instrumented taxis ($N$=60) is examined across Demand Scenarios 1, 2, and 3 for illustrative purposes. Figure \ref{figsensingroutes} depicts all enhanced A* sensing routes generated by vacant instrumented taxis during time period $T$ following the proposed algorithm. Further algorithmic details are provided in Tables \ref{tabmath1} and \ref{tabmath2}. 

	The following observations are derived:

	\begin{itemize}
		\item The origin of each route is defined as the initial position of a vacant instrumented taxi when dispatched for POI monitoring, while the destination is the target  monitoring POI assigned to it. As shown in Figure \ref{figsensingroutes}, the number of distinct destinations is smaller than the total number of monitoring POIs shown in Figure \ref{monitoring POI}. This discrepancy occurs because some monitoring tasks are fulfilled as a byproduct of regular trip assignments to instrumented taxis.
		\item The POI monitoring route is defined as the optimized paths determined by the enhanced A* algorithm to accomplish a sensing task. The distribution of travel demand and the spatial dispersion of monitoring POIs exert a significant influence on these routes of vacant instrumented taxis, which in turn affects the visiting frequencies across different network nodes and edges, as visually summarized in Figure \ref{figsensingroutes}.   
		\item The proposed DRO model quantifies the monitoring cost by applying a unit cost of 1.5 CNY/km to the total distance traveled by vacant instrumented taxis along their assigned POI monitoring routes. This formulation aligns with the comprehensive cost framework, which integrates both the operational expenditures for POI monitoring and the subsidy incentives allocated to the entire instrumented taxi fleet.
	\end{itemize}  
	
	\section{Discussion}\label{secDis}  		
	The current DS ride-hailing market, which comprises mixed instrumented taxi fleets, operates under a conventional deterministic monitoring pattern. This approach leads to highly inefficient operations, primarily due to the excessive monitoring costs incurred by POI monitoring demand.	
	
	The framework proposed in this study serves as a foundational component for a novel mixed DS ride-hailing market design, incorporating an enhanced A* drive-by sensing route policy. Unlike the static and disjointed conventional approach, the proposed framework explicitly accounts for the dynamic fluctuations in POI monitoring, which are driven by various random variables with regard to the market' uncertain ride demand. Furthermore, by optimizing the selection of instrumented taxi fleets according to the heterogeneous spatio-temporal distribution of travel demand, this framework enhances the resource efficiency of the DS ride-hailing market.
	
	Within the paradigm of the enhanced A* drive-by sensing policy, and given the current absence of methods for evaluating dynamic monitoring route effectiveness, even the least sophisticated variant of our framework (the deterministic model) yields a feasible monitoring cost scheme and demonstrates potential for improving upon the status quo. The DRO model and its rolling-horizon variant achieve superior performance against the Design KPI. The distributionally robust technique is particularly effective at managing highly volatile dynamic demand, delivering robust and reliable outcomes. The identified monitoring patterns are validated using ten historical datasets across three distinct demand distributions. Therefore, the validated framework facilitates cost optimization by translating unimpeded POI monitoring demand into quantifiable costs under the enhanced A* drive-by sensing route policy.
	
	The automation of the proposed framework is anticipated to alleviate the operational workload within the DS ride-hailing market, substantially reducing specific tasks such as route planning for vacant instrumented taxis and monitoring cost optimization. Furthermore, the framework's robustness, demonstrated through the DRO approach, allows it to effectively accommodate daily ride demand variability and uncertainty. Although developed and validated for the Longquanyi District, the methodology is readily transferable to other regions with POI monitoring requirements. It should be noted, however, that implementation may shift the operational focus from active sensing to passive sensing, a transition with significant practical implications.	    
	
	\section{Conclusion}\label{secCon} 		
	This paper introduces a DS ride-hailing market that incorporates mixed instrumented taxis to address the demand for tasks monitoring. In contrast to existing approaches, the proposed methodology accounts for the uncertainty arised from uncertain ride demand. It serves the dual purposes of POI monitoring and the robust assessment of monitoring costs. The methodology is summarized as follows:
	\begin{itemize}
		\item An enhanced A* drive-by sensing route service policy is formulated to coordinate vacant instrumented taxis with the monitoring of POIs, within the context of a ride-hailing market operation that utilizes mixed instrumented fleets.
		\item The DRO technique is employed to address ambiguities in the underlying distributions of complex, correlated random variables. This approach simultaneously tackles challenges related to insufficient data availability and limited data granularity. 
	\end{itemize}
	
	To rigorously validate and quantify the efficacy of the proposed framework, extensive air quality monitoring simulations are carried out under a real-world scenario. The results demonstrate that:
	
	\begin{enumerate}
		\item Under the enhanced A* drive-by sensing route service policy, the DRO$_{RH}$ technique demonstrates superior performance compared to the standard Deterministic and DRO techniques. Its advantage lies in prioritizing the reduction of conservatism, as evidenced by the evaluation results of the monitoring cost KPI. For the DRO$_{RH}$ optimization method, the optimal sensor deployment is 60 fleets across all demand scenarios for the specified POIs. Demand Scenario 1, which has a relatively unbalanced spatial trip distribution, exhibits the highest monitoring cost per unit time per vehicle.
		\item The proposed operational framework for the DS ride-hailing platform, integrating optimal sensor deployment and the enhanced A* drive-by sensing route, achieves the anticipated performance in both service quality and sensing capabilities. Service quality is measured by the matching rate and average wait time, while sensing capabilities are evaluated through sensing utility, percentage of grids covered per hour, and monitoring completion rate.  
		\item Planning vacant instrumented taxis to monitor POI makes data collection from monitoring blind spots more timely. The enhanced A*-based drive-by sensing routes between vacant instrumented taxis and monitoring POIs are incorporated into the platform's POI monitoring cost at a rate of 1.5 CNY/km. In addition to these dedicated routes, POI monitoring may also be completed indirectly when instrumented taxis fulfill actual  customer trip requests.
	\end{enumerate}   
	
	Based on travel demand data and targeted POIs, the DRO framework is recommended as an effective guideline for monitoring cost optimization, providing strategic direction for the long-term operation of ride-hailing markets with mixed instrumented fleets. An important extension of the current framework involves incorporating sensor failures and heterogeneous sensor types into the model. Such challenges could be addressed using the DRO methodology presented in this study.	
	 
	\section*{Declaration of competing interest}
	The authors declare that they have no known competing financial interests or personal relationships that could have appeared to influence the work reported in this paper. 
	 	 
	\section*{CRediT authorship contribution statement} 
 	Binzhou Yang: Writing – original draft, Visualization,   Methodology, Investigation, Formal analysis, Conceptualization, Resources. Bin Shuai: Supervision, Resources. Rui Wang: Project administration, Funding acquisition. 
	 
	\section*{Acknowledgement} 
	We thank Dr. Lei Yang for the helpful discussions. 
	This work is partly supported by the Fundamental Research Funds for the Central Universities, China (Grant No. 26CAFUC04040). 

	\section*{Appendix A $\phi$-divergence}\label{Appendix A.} 
	The adjoint and conjugate of $\phi$-divergence for some specific choices.
	\renewcommand\thetable{A.\arabic{table}}
	\setcounter{table}{0}

	\setlength{\LTleft}{0pt}
	\setlength{\LTright}{0pt}

	\begin{longtable}{cccccc}
	\caption{Some specific choices of $\phi$-divergence and their adjoint and conjugate}
	\label{divergencecategory} \\
	
	\hline
	$\phi$-divergence & $\phi(u)$ & $I_{\phi}(p,q)$ & popping & suppression \\
	\hline
	\endfirsthead
	
	\caption*{Table~\thetable{} (continued): Some specific choices of $\phi$-divergence and their adjoint and conjugate} \\
	\hline
	$\phi$-divergence & $\phi(u)$ & $I_{\phi}(p,q)$ & popping & suppression \\
	\hline
	\endhead
	
	\hline
	\multicolumn{5}{r}{{Continued on next page}} \\
	\endfoot
	
	\hline
	\endlastfoot
	
	Kullback-Leibler & $\phi_{kl}(u)=ulogu-u+1$ & $\sum p_{j}log\frac{p_{j}}{q_{j}}$ & $\times$ & $\checkmark$ \\
	Burg entropy & $\phi_{b}(u)=-logu+u-1$ & $\sum q_{j}log\frac{q_{j}}{p_{j}}$ & $\checkmark$ & $\times$ \\
	J-divergence & $\phi_{j}(u)=(u-1)log u$ & $\sum(p_{j}-q_{j})log\frac{p_{j}}{q_{j}}$ & $\times$ & $\times$ \\
	$\chi^{2}$-divergence & $\phi_{c}(u)=\frac{1}{u}(u-1)^{2}$ & $\sum \frac{(p_{j}-q_{j})^{2}}{p_{j}}$ & $\checkmark$ & $\times$ \\
	Modified $\chi^{2}$-divergence & $\phi_{mc}(u)=(u-1)^{2}$ & $\sum \frac{(p_{j}-q_{j})^{2}}{q_{j}}$ & $\times$ & $\checkmark$ \\
	Hellinger distance & $\phi_{h}(u)=(\sqrt{u}-1)^{2}$ & $\sum(\sqrt{p_{i}}-\sqrt{q_{i}})^{2}$ & $\checkmark$ & $\checkmark$ \\
	Variation distance & $\phi_{v}(u)=|u-1|$ & $\sum|p_{i}-q_{i}|$ & $\checkmark$ & $\checkmark$ \\
	\end{longtable}
 
	\section*{Appendix B The enhanced A* drive-by sensing route}\label{Appendix B.}
	
	\renewcommand\thetable{B.\arabic{table}}    
	\setcounter{table}{0}  

	The POI monitoring routings for a single vacant instrumented taxi and for multiple vacant instrumented taxis at dual spatial scales are shown in Table \ref{tabmath1} and Table \ref{tabmath2}, respectively.
		
		\setlength\LTleft{0pt}
	\setlength\LTright{0pt}
	\begin{longtable}{@{\extracolsep{\fill}}rl}
		\caption{Executing a  $N\text{-}\mathrm{route}$   and corresponding $G\text{-}\mathrm{route}$ }  \label{tabmath1}   \\
		\hline
		\multicolumn{2}{l}{Algorithm1}     \\ 
		\hline
		\\ 
		Input   &  Starting node $\tau_{1}\in g_{1}$, corresponding monitoring POI $\tau\in\mathcal{N}$.
		\\
		Step one    &  \begin{tabular}[c]{@{}l@{}} Apply the modified A* algorithm to find  the shortest path $p^{*}=(\tau_{1},...,\tau_{m})$ and \\ corresponding $N\text{-}\mathrm{route}$   $\bar{\kappa}(\tau_{1};\kappa_{G})$.   \end{tabular} 
		\\
		Step two & Calculate the actual $G\text{-}\mathrm{route}$  $G_{\bar{\kappa}(\tau_{1};\kappa_{G})}$ traversed by $\bar{\kappa}(\tau_{1};\kappa_{G})$.
		\\
		Step three & Calculate $\varphi(a\tau)$ and matching pair $x_{a\tau}$ combined with Eqs. (\ref{ueffi}) - (\ref{ueffiopt4}) .
		\\
		Output & Matching pair $x_{a\tau}$, $N\text{-}\mathrm{route}$  $\bar{\kappa}(\tau_{1};\kappa_{G})$, $G\text{-}\mathrm{route}$   $G_{\bar{\kappa}(\tau_{1};\kappa_{G})}$.
		\\
		\\\hline
	\end{longtable} 

	    \setlength{\LTleft}{0pt}
	\setlength{\LTright}{0pt}
	
	\begin{longtable}{@{\extracolsep{\fill}}rl}
		\caption{Multi-vehicle routing for discrete time window $t$}
		\label{tabmath2} \\
		
		\hline
		\multicolumn{2}{l}{Algorithm 2} \\
		\hline
		\endfirsthead
		
		\caption*{Table~\thetable{} (continued): Multi-vehicle routing for discrete time window $t$} \\
		\hline
		\multicolumn{2}{l}{Algorithm 2 (continued)} \\
		\hline
		\endhead
		
		\hline
		\multicolumn{2}{r}{{Continued on next page}} \\
		\endfoot
		
		\hline
		\endlastfoot
		\\
		Input   &\begin{tabular}[c]{@{}l@{}}  Vacant instrumented taxis $n=1...n^{A}$ and monitoring POI $\tau\in\mathcal{N}$, \\ maximum iteration limit max\underline{ }iter. \end{tabular}
		\\
		Initialize   & Set $k=1$.
		\\
		Step one    &  \begin{tabular}[c]{@{}l@{}} Apply Algorithm 1 to the vacant instrumented taxis $n=1...n^{A}$ and \\ monitoring POI $\tau\in\mathcal{N}$.  \end{tabular}
		\\
		Step two & \begin{tabular}[c]{@{}l@{}}  Adjustment of routing solutions.\\For $n=1,...,n^{A}$: \end{tabular}
		\\
		Step two (a) & \begin{tabular}[c]{@{}l@{}}  \quad\quad  If the
			$n$-th vacant instrumented taxi has formed a monitoring POI \\\quad\quad matching pair, the pair is dissolved and the coverage utility contributed  \\\quad\quad by the
			$n$-th vacant instrumented taxi is removed; \end{tabular}
		\\
		Step two (b) & \begin{tabular}[c]{@{}l@{}}  \quad\quad Apply the Algorithm 1 to the $n$-th vacant instrumented taxi, if the new  \\ \quad\quad solution  yields higher  utilities than the current solution, then accept the\\ \quad\quad new solution; otherwise, retain the current solution; set $k$=$k$+1;
		\end{tabular}
		\\
		Step three & \begin{tabular}[c]{@{}l@{}}  \quad\quad If no solutions are updated for any of the vacant instrumented taxis or\\ \quad\quad $k$=max\underline{ }iter, terminate the Algorithm 2; otherwise, go to Step two.\end{tabular}
		\\
		Output & Routing solutions $R^{A}={R_{1}^{A},...,R_{n^{A}}^{A}}$ for $n^{A}$ vacant
		instrumented taxis.
		\\
		\\
	\end{longtable}
	
	\section*{Appendix C Operational performance figures}\label{Appendix C.}

	In line with Section \ref{4.3}, the operational performance across three demand scenarios for fleet sizes of 10-60 vehicles is presented in Figures \ref{figsim3casesm}--\ref{figsim3casess}, including level of service, sensing externality and monitoring capability. In particular, building upon prior research \citep{YHLL2025}, we have incorporated simulation tests of monitoring capability in Figure \ref{figsim3casess}.
	
	\renewcommand\thefigure{C.\arabic{figure}}    
	\setcounter{figure}{0}  
	
	\begin{figure}[!h]
		\centering
		\includegraphics[scale=0.48,height=0.22\textheight]{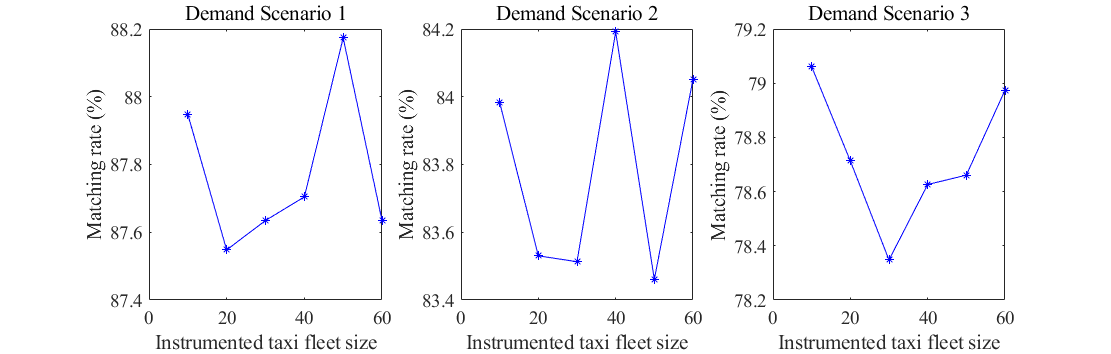}\\
		\caption{Matching rate under different demand scenarios and instrumented taxi fleet size.}\label{figsim3casesm} 
	\end{figure}

	\begin{figure}[!h]
		\centering
		\includegraphics[scale=0.50,height=0.22\textheight]{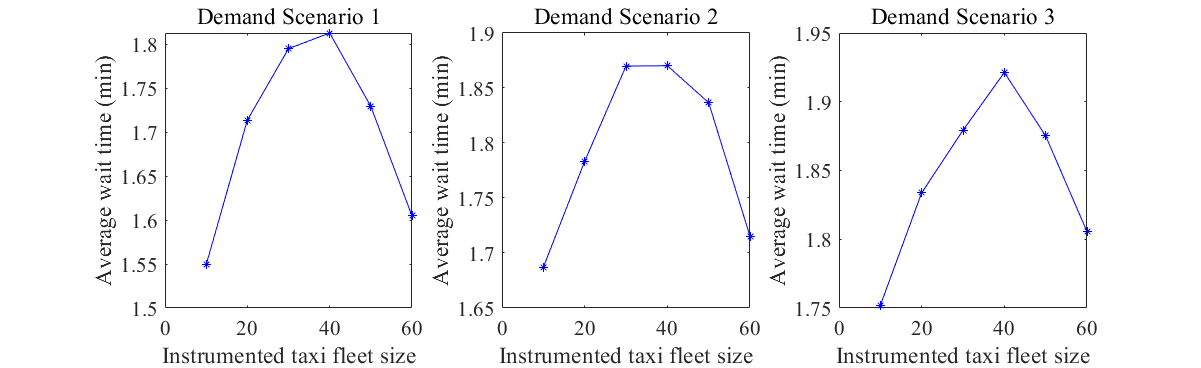}\\
		\caption{Average wait time under different demand scenarios and instrumented taxi fleet size.}\label{figsim3casesw} 
	\end{figure} 

	\begin{figure}[!h]
		\centering
		\includegraphics[scale=0.55,height=0.40\textheight]{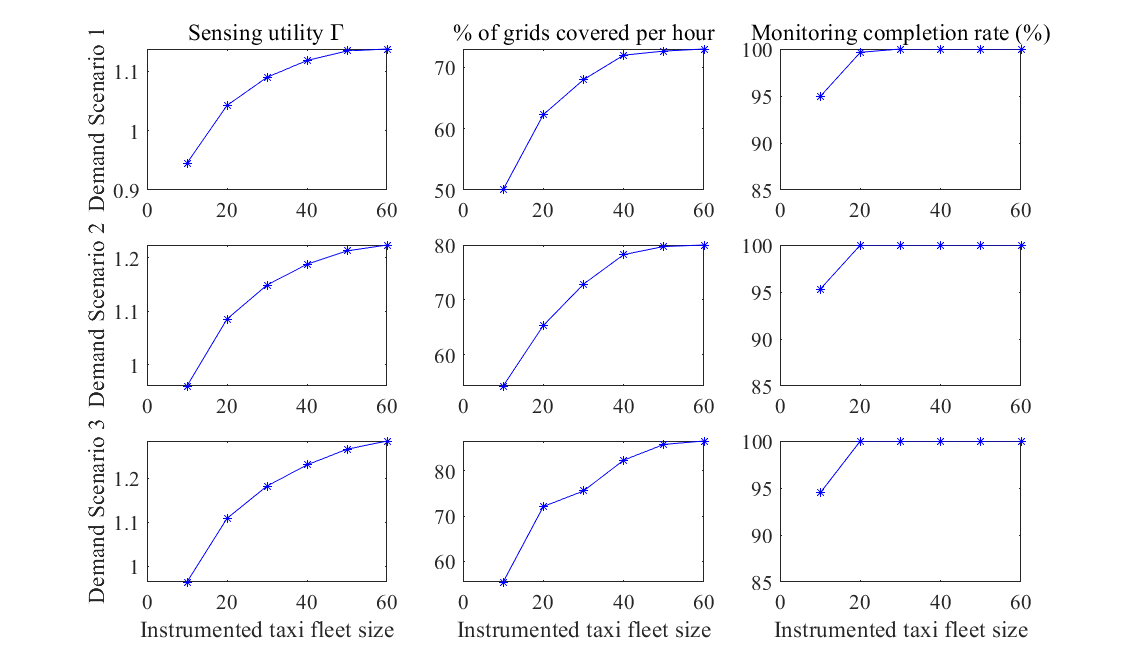}\\
		\caption{The sensing externality and monitoring capacity under different demand scenarios and instrumented taxi fleet size.}\label{figsim3casess} 
	\end{figure}

    \pagebreak
	
\end{document}